%% file: main.tex
\documentclass{amsart}
\pdfoutput=1
\usepackage[utf8]{inputenc}
\usepackage[english]{babel}
\usepackage{amsmath,amssymb,amsthm,graphicx,url,color,enumerate,dsfont,stmaryrd}
\usepackage{xpatch}
\usepackage{IEEEtrantools}
\usepackage{mathtools}
\usepackage{tikz}
\usetikzlibrary{calc, arrows.meta, patterns, positioning, decorations, decorations.markings}
\usepackage{subcaption}
\usepackage[colorlinks=true,citecolor=blue]{hyperref}
\usepackage{cleveref}
\usepackage[top=4cm,bottom=4cm,margin=3cm]{geometry}
\usepackage{csquotes}
\usepackage[protrusion=true,expansion,babel=true]{microtype}

\theoremstyle{plain}
\newtheorem{theorem}{{Theorem}}[section] 
\newtheorem*{theorem*}{{Theorem}}
\newtheorem{proposition}[theorem]{Proposition}
\newtheorem*{proposition*}{Proposition}
\newtheorem{corollary}[theorem]{Corollary}
\newtheorem*{corollary*}{Corollary}
\newtheorem{lemma}[theorem]{Lemma}
\newtheorem*{lemma*}{Lemma}

\theoremstyle{definition}
\newtheorem{definition}[theorem]{Definition}
\newtheorem*{definition*}{Definition}

\theoremstyle{remark}
\newtheorem{remark}[theorem]{Remark}

\makeatletter

\@addtoreset{equation}{section}
\makeatother

\newcommand{\abs}[1]{\left\vert #1\right\vert}      
\newcommand{\nr}[1]{\left\Vert #1\right\Vert}          
\newcommand{\innp}[2]{\left \langle #1 , #2 \right \rangle}           

\newcommand{\set}[1]{\left\{ #1 \right\}}		

\renewcommand{\leq}{\leqslant}	\renewcommand{\geq}{\geqslant}
\newcommand{\inv}{^{-1}}

\newcommand{\littleo}[3][]{\mathop{o}\limits_{\ifblank{#1}{#2 \to #3}{\substack{#2\to #3\\#1}}}}
\newcommand{\1}{\mathds 1}       
\let\Re\relax \DeclareMathOperator{\Re}{Re}        
 
\newcommand{\limt}[3][]{\xrightarrow[\ifblank{#1}{#2 \to #3}{\substack{#2\to #3\\#1}}]{}}

\DeclareMathOperator{\Dom}{Dom}
\DeclareMathOperator{\Rank}{Rank}
\DeclareMathOperator{\Sp}{\s}
\DeclareMathOperator{\Span}{span}

\DeclareMathOperator{\distance}{distance}

\DeclareMathOperator{\supp}{supp}
\DeclareMathOperator{\Ran}{Ran}

\newcommand{\R}{\mathbb{R}}		\newcommand{\C}{\mathbb{C}}
\newcommand{\N}{\mathbb{N}}	\newcommand{\Z}{\mathbb{Z}}	  
\newcommand{\T}{\mathbb{T}}

\newcommand{\Lc}{{\mathcal L}}
\newcommand{\Dc}{{\mathcal D}}
\newcommand{\Sc}{{\mathcal S}}
\newcommand{\Zc}{\mathcal Z}
\DeclareMathOperator{\Hol}{\mathcal O}

\renewcommand{\a}{\alpha}\renewcommand{\b}{\beta}\newcommand{\g}{\gamma}\renewcommand{\d}{\delta}\newcommand{\e}{\varepsilon}\newcommand{\z}{\zeta}\renewcommand{\th}{\theta}\newcommand{\Th}{\Theta}\renewcommand{\l}{\lambda}\newcommand{\s}{\sigma}\newcommand{\f}{\varphi}\renewcommand{\o}{\omega}\renewcommand{\O}{\Omega}

\mathtoolsset{mathic}
\newcommand*{\noic}{\sb{}\kern-\scriptspace }

\makeatletter
\patchcmd{\@IEEEeqnarray}{\relax}{\relax\intertext@}{}{}
\makeatother

\newcounter{proofstep}[theorem]
\newcommand{\step}[1]{%
\refstepcounter{proofstep}%
\vskip-\lastskip\medskip\noindent\textit{Step \arabic{proofstep}: #1.}}
\xapptocmd\proof{\setcounter{proofstep}{0}}{}{}

\newcommand{\diff}[1][-3]{\ensuremath{\mathop{}\mkern#1mu\mathrm{d}}}
\newcommand{\eu}{\mathrm e}
\newcommand{\iu}{\mathrm i}

\DeclareMathOperator{\dagmon}{d_{\mathrm{Agm}}}
\DeclareMathOperator{\dagmontilde}{\delta}

\newcommand{\He}{h} %Hermite polynomials
\newcommand{\sector}[1][\theta_0]{\Sigma_{#1}}
\newcommand{\pol}{{\mathrm{pol}}}
\newcommand{\h}{{\mathrm{H}}}
\newcommand{\Op}[1]{\mathop{#1(z\partial_z)}}

\DeclareMathOperator{\Tr}{Tr}

\begin{document}

%TODO: 
% -Dessin annexe cutoff pour chemin simple
% -???
% -Profit

\title[Null-controllability of the Baouendi-Grushin equation]{Null-controllability properties of the generalized two-dimensional Baouendi-Grushin equation with non-rectangular control sets}

\author{Jérémi Dardé}
\author{Armand Koenig}
\author{Julien Royer}

\address{Institut de Math\'ematiques de Toulouse, UMR 5219, Universit\'e de Toulouse, CNRS,  UPS, F-31062 Toulouse Cedex 9, France.}

\email[J. Dardé]{jeremi.darde@math.univ-toulouse.fr}
\email[A. Koenig]{armand.koenig@math.univ-toulouse.fr}
\email[J. Royer]{julien.royer@math.univ-toulouse.fr}

\subjclass[2010]{35K65, 93B05, 47B28, 47A10}

\begin{abstract}
We consider the null-controllability problem for the \emph{generalized Baouendi-Grushin equation} $(\partial_t - \partial_x^2 - q(x)^2\partial_y^2)f = \1_\o u$ on a rectangular domain. Sharp controllability results already exist when the control domain $\o$ is a vertical strip, or when $q(x) = x$. In this article, we provide upper and lower bounds for the minimal time of null-controllability for general $q$ and non-rectangular control region $\o$. In some geometries for $\o$, the upper bound and the lower bound are equal, in which case, we know the exact value of the minimal time of null-controllability.

Our proof relies on several tools: known results when $\o$ is a vertical strip and cutoff arguments for the upper bound of the minimal time of null-controllability; spectral analysis of the Schrödinger operator $-\partial_x^2 + \nu^2 q(x)^2$ when $\Re(\nu)>0$, pseudo-differential-type operators on polynomials and Runge's theorem for the lower bound.
\end{abstract}

\maketitle

\section{Introduction and statements of the main results}

\subsection{The Baouendi-Grushin equation}

In this article, we study some controllability properties of the two-dimensional \emph{generalized Baouendi-Grushin equation}. 

Let $L_-,L_+$ be positive, and $I = (-L_-,L_+)$. Let $q \in C^0(\overline I)$ such that $q(0)= 0$ and $q(x) \neq 0$ for all $x \in \overline I\setminus \{0\}$. We denote by
$\T$ the one-dimensional torus $\R/2\pi \Z$. For $T > 0$, $f_0 \in L^2(I\times \T)$ and $F \in L^2((0,T);L^2(I \times \T))$ we consider the generalized Baouendi-Grushin equation

\begin{equation}
 \label{eq-grushin-gen-intro}
 \left\{\begin{array}{ll}
         (\partial_t-\partial_x^2 - q(x)^2 \partial_y^2)f(t,x,y) = F(t,x,y),&\quad t\in (0,T), x\in I, y\in \T\\
         f(t,x,y) = 0,&\quad t\in (0,T), x \in \partial I, y\in \T, \\
         f(0,x,y) = f_0, & \quad x\in I, y\in \T.
        \end{array}\right.
\end{equation}

Note that, because $q(0) = 0$, the equation degenerates on the vertical axis $\left\lbrace 0 \right\rbrace \times \T$. Nevertheless, the equation is well posed.
Precisely, the Friedrichs extension (see \cite[Section 4.3]{Helffer13})  of the operator
\[
 f \in C^\infty_c(I\times \T) \mapsto - \partial^2_x f - q(x)^2 \partial_y^2 f,
\]
generates an analytic semigroup, which allows to define a solution 
of the generalized Baouendi-Grushin equation~\eqref{eq-grushin-gen-intro} in the sense of semigroups \cite{Brezis11}. In our case, this solution is smooth in the following sense:

\begin{proposition} \label{prop_wellposedness}
	For any source term $F \in L^2((0,T);{L^2(I\times \T)})$ and any initial condition $f_0 \in L^2(I\times \T)$, there exists a unique $f \in C^0([0,T];L^2(I\times \T)) \cap L^2((0,T);V)$
	solution of the generalized Baouendi-Grushin equation~\eqref{eq-grushin-gen-intro}, with 
	\[
	V = \left\lbrace f \in L^2(I\times \T), \ \partial_x f \in L^2(I\times \T),\ q\, \partial_y f \in L^2(I\times \T) \right\rbrace.
	\] 
\end{proposition}

This result is proved in \cite{BCG14} in the case $q(x) = |x|^\gamma$ with $\gamma>0$. The proof is easily generalized to our case of interest.

\subsection{Control problem for the Baouendi-Grushin equation}

Our study focuses on internal null-controllability of the Baouendi-Grushin equation. More
precisely, let $\omega \subset I \times \T$ be a non-empty open set and $u \in L^2((0,T); \, L^2(\o))$. 
The controlled Baouendi-Grushin equation reads
\begin{equation}
\label{eq-grushin-gen-contr}
\left\{\begin{array}{ll}
(\partial_t-\partial_x^2 - q(x)^2 \partial_y^2)f(t,x,y) = \mathds 1_\o u(t,x,y),&\quad t\in (0,T), x\in I, y\in \T\\
f(t,x,y) = 0,&\quad t\in (0,T), x \in \partial I, y\in \T, \\
f(0,x,y) = f_0, & \quad x\in I, y\in \T,
\end{array}\right.
\end{equation}
where $f$ is the state of the system, and $\1_\o u$ is the control supported in $\o$. 

\begin{definition}[Null-controllability] Let $T>0$. The Baouendi-Grushin equation \eqref{eq-grushin-gen-contr} is said to be null-controllable on $\o$ in time $T$ if, for any initial condition 
	$f_0 \in L^2(I\times \T)$, there exists $u \in L^2((0,T);L^2(\o))$ such that the solution $f$ of \cref{eq-grushin-gen-contr} satisfies $f(T,\cdot,\cdot) = 0$ in $I\times \T$.
\end{definition}

It is known that, contrary to usual non-degenerate parabolic equations like the heat equation, due to the degeneracy of $q$ on $\left\lbrace 0 \right\rbrace \times \T$, the null-controllability properties of \eqref{eq-grushin-gen-contr} strongly depend on the control set $\o$ and the time horizon $T$. More precisely, for certain control sets $\o$, there is no time $T>0$ such that \cref{eq-grushin-gen-contr} is  null-controllable, whereas for other control sets $\o$ a minimal time of null-controllability appears. We refer to the bibliographical comments, \cref{section_biblio} below, for a detailed 
description of the known results on the subject.

In the present paper, we aim to give precise null-controllability results for equation \eqref{eq-grushin-gen-contr}, for a large class of control sets $\o$ and a large class of functions $q$.

\subsection{Main results}
We are interested in the case where the equation is degenerate on $\{x=0\}$. Thus, we assume that $q(x) = 0$ only when $x=0$, and we assume without loss of generality that $q'(0) > 0$.

Before presenting the main results of our study, we introduce 
the so-called \emph{Agmon distance} of a point $x \in I$
to the origin, defined by:\footnote{It corresponds to the usual definition of the Agmon distance given for example in~\cite[Eq.~(6.3)]{DS99}, with $V = q^2$ and $E = 0$.}
\begin{equation} \label{eq-def-agmon}
\dagmon : x\in I  \mapsto \int_0^x q(s)\diff s.
\end{equation}
This quantity appears naturally in the computation of the minimal time of null-controllability for the generalized Baouendi-Grushin equation.

\subsubsection{Lack of null-controllability in small time for a class of control sets $\o$}

Our main result is a negative result of null-controllability for small times. We show that if the control set $\o$ stays at positive distance from a horizontal segment of the form $(a,b)\times \left\lbrace y_0 \right\rbrace$, with
$-L_-\leq a<0<b\leq L_+$ and $y_0 \in \T$, then equation \eqref{eq-grushin-gen-contr} is not null-controllable on $\o$ for time $T$ smaller than a precisely given critical time.

To properly state the result, we need to introduce a modified version of the Agmon distance. We set $\dagmontilde(x) = \dagmon(x)$ for $x \in I$, and
$\dagmontilde(-L_-) = \dagmontilde(L_+) = +\infty$.

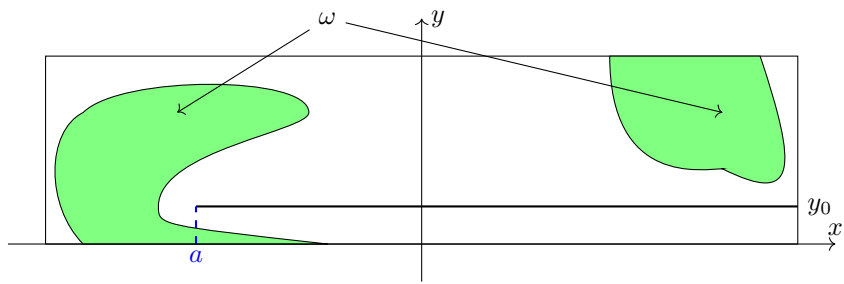
\begin{figure}[h]
\centering
\input{omega_neg_intro}
\caption{In green, an example of a domain $\omega$ with, in thick black, an example of a horizontal segment that stays at positive distance of $\o$. In this example, \cref{th-main-neg-intro} implies that the generalised Baouendi-Grushin equation is not null-controllable in time $T<\dagmon(a)$.}
\label{fig:omega_neg_intro}
\end{figure}

\begin{theorem} \label{th-main-neg-intro}
	Assume that $q\in C^2(\overline I)$ is such that $q(0) = 0$, $q'(0) > 0$ and $q(x) \neq 0$ whenever $x\neq 0$.
	Let $\o$ be an open subset of $I\times \T$. Assume that there exist $a \in [-L_-,0)$, $b \in (0,L_+]$ and $y_0 \in \T$ such that 
	\[
	\distance\big((a,b)\times \{y_0\}, \o\big) > 0.
	\]
	(See \cref{fig:omega_neg_intro}.) Then, the generalized Baouendi-Grushin equation~\eqref{eq-grushin-gen-contr} is not null-controllable on $\o$ in time $T$ such that 
	\[
	T< \frac{1}{q'(0)}\min \left(\dagmontilde(a),\dagmontilde(b) \right).
	\]
\end{theorem}

This theorem is a generalization of \cite[theorem 3.3]{DK20}, where the result is proved for $I$ symmetric with respect to the origin and, more restrictively, for $q(x) = x$.  
A key step in our proof of \cref{th-main-neg-intro} is the study of spectral properties of the family of operators defined on $L^2(I)$ by
\begin{equation} \label{eq_def_Pnu}
P_\nu\colon -\partial^2_x + \nu^2 q(x)^2, \quad \Dom(P_\nu) = H^2(I)\cap H^1_0(I).
\end{equation}
Because of technical reasons, we have to consider $P_\nu$ for every $\Re(\nu)>0$, which makes $P_\nu$ non self-adjoint.

In the previous article~\cite{DK20}, the corresponding results were obtained using explicit solutions of particular ordinary differential equations. These explicit formulae are not available in our general setting.

In our case, we obtain a localization of an eigenvalue of $P_\nu$, as well as precise Agmon type estimates for an associated eigenfunction,\footnote{For $\nu \in \R_+^*$, we localize the smallest eigenvalue of $P_\nu$, and for $\nu$ complex, we localize its analytic continuation.} uniformly in $\nu = \vert \nu \vert \eu^{\iu \theta}$ with $\vert \nu \vert$ large enough and $0 \leq \vert \theta \vert  \leq \theta_0$ for some $\th_0 \in [0, {\pi}/{2} )$. To that end, we compare $P_\nu$ with the non-selfadjoint harmonic oscillator $H_{q'(0)\nu}\coloneqq -\partial_x^2 + (q'(0)\nu)^2 x^2$.
\Cref{th-main-neg-intro} is proved in \cref{sec_negative_result}, with the spectral analysis done in \cref{sec-spectral}.

Note that \cite{HSS05} contains closely related spectral asymptotics. However, we cannot apply them directly since our domain has a boundary and we need uniform estimates with respect to the parameter $\nu$.

\subsubsection{Precise critical time of null-controllability for a class of control sets $\o$}
With \cref{th-main-neg-intro} we can deduce in some setting the critical time for the null-controllability of \eqref{eq-grushin-gen-contr}.

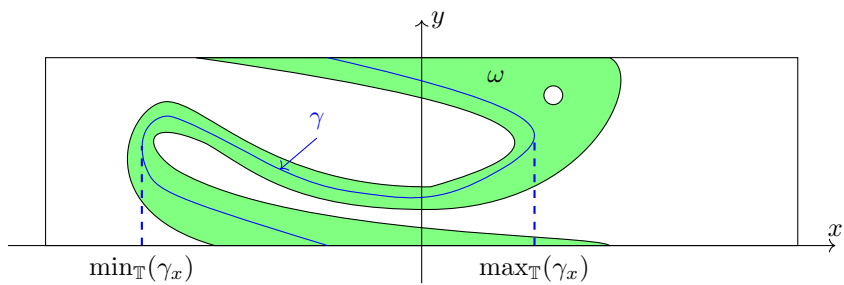
\begin{figure}[h]
\centering
\input{omega-pos-intro}
\caption{In green, an example of a domain $\o$, with, in blue, a corresponding path $\g$ that satisfies the hypotheses of \cref{th-main-pos-intro}.}
\label{fig:omega-pos-intro}
\end{figure}
We first mention a natural adaptation of \cite[theorem 3.1]{DK20}. It was actually claimed in~\cite[Remark 3.2]{DK20}, but the statement was imprecise if $q$ is not odd. We take the opportunity to correct the statement:
\begin{theorem} \label{th-main-pos-intro}
	Assume that $q\in C^3(\overline I)$ is such that $q(0) = 0$ and $\min_{\overline I} q' >0$. Let $\o$ be an open subset of $I\times \T$. Assume that there exists a closed path $\gamma = (\g_x, \g_y) \in C^0(\T;\, \omega)$ 
	such that $\{-L_-\}\times \T$ and $\{L_+\}\times \T$ are included in different connected components of $(\overline I\times\T) \setminus \g(\T)$ (see \cref{fig:omega-pos-intro}).
	
	The generalized Baouendi-Grushin equation~\eqref{eq-grushin-gen-contr} is null-controllable on $\o$ in time $T$ such that 
	\[
	T> \frac{1}{q'(0)}\max\left(\dagmon\left(\min_\T(\g_x)\right), \dagmon\left( \max_\T (\g_x)\right)\right).
	\]
\end{theorem}
\begin{remark}\label{rk-homotopy}
In this theorem, we can replace the hypothesis ``$\{-L_-\}\times \T$ and $\{L_+\}\times \T$ are included in different connected components of $(\overline I\times\T) \setminus \g(\T)$'' by ``$\g$ is not homotopic to a constant path''. These two conditions are essentially equivalent. We discuss this in \cref{th-component,th-homotopy-2,rk-homotopy-2}.
\end{remark}

We detail the proof of this theorem in \cref{sec_positive_result}.
Combining \cref{th-main-neg-intro,th-main-pos-intro}, we obtain the following result which gives the precise critical time of null-controllability of the generalized Baouendi-Grushin
for a large class of control sets $\o$ and functions $q$:

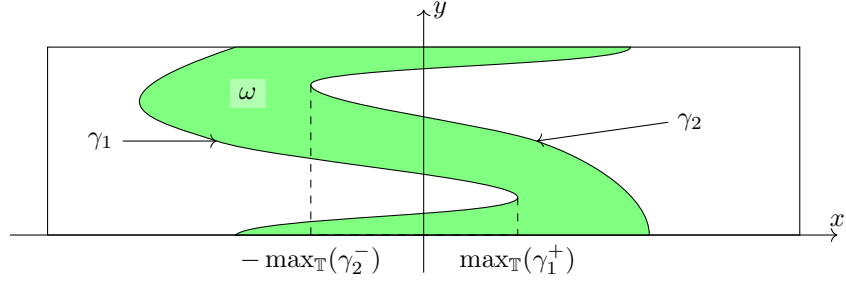
\begin{figure}
\centering
\input{omega-corollary-intro}
\caption{In green, an example of a domain $\omega$ that satisfies the hypotheses of \cref{th-main}.}
\label{fig:omega_corollary_intro}
\end{figure}
\begin{theorem}\label{th-main}
 Assume that $q\in C^3(\overline I)$ is such that $q(0) = 0$ and $\min_{\overline I} q' > 0$. Let $\g_1,\g_2\in C^0(\T;\, I)$ such that for every $y\in \T$, $\g_1(y) < \g_2(y)$. Let $\omega = \{(x,y)\in I\times \T : \g_1(y)< x < \g_2(x) \}$ (see \cref{fig:omega_corollary_intro}) and\footnote{When $f$ is a real valued function, we denote $f^+ = \max(f,0)$ and $f^- = \max(-f,0)$ its positive and negative part respectively.}
 \[
  T_* \coloneqq \frac{1}{q'(0)} \max\left(\dagmon\left(-\max_\T(\g_2^-)\right),\dagmon\left(\max_\T(\g_1^+)\right)\right).
 \]
Then the generalized Baouendi-Grushin equation~\eqref{eq-grushin-gen-contr} is null-controllable on $\omega$ in any time $T > T_*$, but it is not null-controllable on $\o$ in time $T<T_*$.
\end{theorem}
This theorem is proved in \cref{sec-crit-time}.

\subsubsection{Comments}

Before proceeding further, we make some additional comments on our results.

\begin{itemize}
	\item The assumptions regarding the function $q$ in \cref{th-main-neg-intro} are slightly more general than in \cref{th-main-pos-intro}. They seem also more natural in the context of our study. Therefore, we conjecture that  \cref{th-main-pos-intro} holds for functions $q \in C^3(\overline I)$ satisfying
	\[
	q(0) = 0, \ q'(0) > 0, \ q(x) \neq 0 \text{ for all } x \in \overline{I}.
	\]
	But up to our knowledge, this is still an open question.
	\item There are still numerous geometrical configurations not included in  \cref{th-main}. Nevertheless, in many situations,  \cref{th-main-pos-intro} and \cref{th-main-neg-intro}  give information about null-controllability properties. As an example, in the geometrical configuration described in  \cref{fig:example3}, combining  \cref{th-main-pos-intro} and  \cref{th-main-neg-intro}, we obtain the existence of a critical time
	\[
	T_* \in \left( \frac {\dagmon(-a)}{q'(0)}, \frac {\dagmon(-b)}{q'(0)} \right)
	\]
	such that the Baouendi-Grushin equation is null-controllable on $\o$ in time $T>T_*$, and is not null-controllable on $\o $ in time $T< T_*$. 
		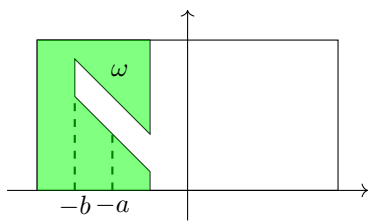
\begin{figure}[h]
		\centering
		\input{example3.tex}
		\caption{In this configuration, we obtain lower and upper bounds of the critical time of null-controllability}
		\label{fig:example3}
	\end{figure}
	In the two geometrical configurations presented in \cref{fig:examples1and2},  \cref{th-main-neg-intro} implies that the Baouendi-Grushin equation is not
	null-controllable on $\o$ in time $T< \dagmon(a)/q'(0)$. Note that in these two configurations, the question of the null-controllability of \eqref{eq-grushin-gen-contr} in $\o$ for some time $T$ large enough is still open.
	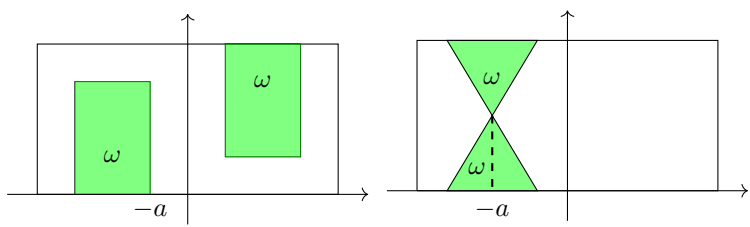
\begin{figure}[h]
		\centering
		\input{example1.tex} \input{example2.tex}
		\caption{In these configurations, we obtain a lower bound on the critical time of null-controllability.}
		\label{fig:examples1and2}
	\end{figure}
	\item These results are stated for the generalized Baouendi-Grushin equation posed on $I\times \T$. They can be adapted to the equation posed on $I\times (0,\pi)$ with Dirichlet boundary conditions, with very similar proofs. We refer to \cref{app-dirichlet} for details on the statements and the corresponding proofs.
\end{itemize}

\subsection{Bibliographical comments} \label{section_biblio}

\subsubsection{On the Baouendi-Grushin equation}

The study of controllability properties of system \eqref{eq-grushin-gen-contr} began with the pioneering work \cite{BCG14}, where the authors study the null-controllability of the equation
\begin{equation} \label{eq_Grushin_x}
\left\{\begin{array}{ll}
(\partial_t  - \partial^2_{x} - \vert x\vert^{2\gamma} \partial^2_y) f 
= 1_\omega u, & t\in (0,T), \ x\in (-1,1),\ y \in \T,\\
 f(t,x,y) = 0, & t \in (0,T), \ x = \pm 1, \ y \in \T. 
\end{array}\right.
\end{equation}
They prove that in the case $\gamma \in (0,1)$ (weak degeneracy), the Baouendi-Grushin equation~\eqref{eq-grushin-gen-contr} is null-controllable for any control set $\o$ and any time $T>0$, whereas in the case $\gamma>1$ (strong degeneracy), it is not null-controllable for any control set $\omega$ and any time $T>0$, except if $\o$ contains $\left\lbrace 0 \right\rbrace \times \T$ in which case it is null-controllable in any positive time $T$.

More surprisingly, in the case $\gamma = 1$, which corresponds to the Baouendi-Grushin equation~\eqref{eq-grushin-gen-contr} with $q(x) = x$, and for  $\omega = (a,b) \times \T$, with $0<a<b$, there exists a critical time
$T_* \geq \frac{a^2}{2}$ such that the Baouendi-Grushin equation~\eqref{eq-grushin-gen-contr}
is null controllable on $\o$ in time $T$, for every $T> T_*$,
and is not null controllable on $\o$ in time $T$, for every $T < T_*$.
It is also proved that if $\gamma = 1$ and $\o$ contains the vertical line $\left\lbrace 0 \right\rbrace \times \T $, equation \eqref{eq_Grushin_x} is null controllable in any time $T>0$. Such a minimal time of null-controllability would not be surprising for equations with finite speed of propagation, such as the wave equation~\cite{BLR92}, but the Baouendi-Grushin equation has a \emph{infinite} speed of propagation.

Many works followed, trying to characterize precisely the critical time $T_*$, and to generalize the result to  different geometrical settings and different functions $q$. The first exact characterization of $T_*$ is given in \cite{BMM15} in the case $q(x) = x$ and with two symmetric vertical strips as control set, that is $\o = (-1, - a) \times (a,1)$, $a \in (0,1)$. Using the transmutation method and sideways energy estimates, the authors prove that \cref{eq-grushin-gen-contr} is null-controllable in $\o$ in any time $T> \frac{a^2}{2}$, and is not null-controllable in $\o$ in any time $T<\frac{a^2}{2}$. 
 
When $\omega$ is a vertical strip of the form $(a,b) \times \T$, with $a>0$, as in \cite{BCG14}, the precise value of the critical time $T_*$ was obtained independently in the works \cite{ABM21,BDE20,LL22a}. More precisely, in \cite{ABM21}, using new estimates for biorthogonal sequences to real exponentials
and the moments method, the authors prove  that in the case $q(x) = x$, the critical time  is $\frac{a^2}{2}$. In \cite{BDE20}, with a function $q$ satisfying the assumptions of  \cref{th-main}, the authors use a Carleman strategy to obtain that \cref{eq-grushin-gen-contr} is null controllable on $\o$ in any time $T> T_*$, and not null-controllable on $\o$ in any time $T < T_*$, with
\[
 T_* = \frac{\dagmon(a)}{q'(0)}.
\]
Very recently, this result was obtained in \cite{LL22a} using the moments method, with a stronger smoothness assumption on $q$ (see
 \cite[Remark 1.12 and Proposition 1.13]{LL22a}).
 
All the strategies developed in \cite{BCG14,BMM15,ABM21,BDE20}, although very different, rely on a Fourier expansion of system \eqref{eq-grushin-gen-contr} with respect to the $y$-variable and the study of the obtained family of one dimensional parabolic equations in the variables $t,x$. As a consequence, the control set $\omega$ has to contain a vertical strip, which seems to be  an important restriction of the proposed methods. 
Nevertheless, in \cite{DK20}, the authors generalize the positive null-controllability results obtained in \cite{BDE20} to a large class of control sets: 
in the setting of \cref{th-main} and with the additional assumptions that $I$ is symmetric and $q$ is odd, system \eqref{eq-grushin-gen-contr} is null controllable in any time $T>T_*$, with 
\[
 T_* = \dfrac1{q'(0)}\max\left(\dagmon(-\max(\g_2^-)),\dagmon(\max(\g_1^+))\right).
\]
In the specific case $q(x) = x$, they also prove that if there exist $a,b \in I$, $a < 0 <b$, and $y_0 \in \T$ such that  $\distance((a,b) \times \left\lbrace y_0 \right\rbrace \cap \omega)  >  0$, then \eqref{eq-grushin-gen-contr} is not null controllable in time $T< \min(a^2, b^2)/2$, whereas if there is $y_0 \in \T$ such that  $\distance (I\times \lbrace y_0 \rbrace \cap \o) > 0$, then \eqref{eq-grushin-gen-contr}  is not null controllable on $\o$ in any positive time $T$. 
 \Cref{th-main} is the generalization of this result to a wider class of functions $q$. 
 
To end this overview on controllability issues for the parabolic Baouendi-Grushin equation, we point out that partial controllability results are known in some multidimensional configurations~\cite{BDE20} while precise results are known for cascade systems of two-dimensional Baouendi-Grushin equations with one control, in the case $q(x) = x$ \cite{ABM21}.

\subsubsection{Some related problems}

Let us briefly mention the literature on related problems, in several directions: other degenerate parabolic equations, minimal time of null controllability for parabolic systems, and other type of degenerate equations.

Since the pioneering works \cite{Egorov63,FR71} on the null-controllability of the one-dimensional heat equation, the null-controllability of non-degenerate parabolic equations has been extensively studied. The null-controllability of degenerate
parabolic equations is a more recent subject of study. The case of a degeneracy at the boundary of the domain is now well-understood~\cite{CMV16} (see also the references therein). 

When the degeneracy occurs in the domain, we lack for the moment a general theory, and equations are studied case by case. The two-dimensional Baouendi-Grushin equations is arguably the simplest and best understood equation of that type.
Very similar results, including a minimal time of null-controllability for quadratic degeneracy, have been observed for the heat equation on the Heisenberg group \cite{BC17,BDE20}, and the Kolmogorov equation \cite{Beauchard14,BZ09,BHHR15,DR21,Koenig20}. 

The related problem of \emph{approximate controllability} for degenerate parabolic equations has been studied in a somewhat general framework~\cite{LL22}.

A minimal time of null-controllability might also appear for the heat equation with punctual control~\cite{Dolecki73} and for systems of parabolic equations, degenerate or not \cite{ABGd16,BBM20}.

Finally, let us mention than the subelliptic wave equation is not controllable~\cite{Letrouit20}, and that the Grushin-Schrödinger equation has a minimal time of controllability~\cite{BS19,LS20}.

\section{Null-controllability in large time} \label{sec_positive_result}

In this section, we prove \cref{th-main-pos-intro}.
The idea of the proof is to use known controllability results for equation \eqref{eq-grushin-gen-contr} when the control set is a vertical strip combined with a cutoff argument. More precisely, we recall the following result ~\cite[theorem 1.4]{BDE20}.\footnote{The reference \cite{BDE20} states the result with a control on the boundary. But cutoff arguments allow to construct controls on vertical strips from boundary controls, as in \cite[Appendix A]{DK20}.}
\begin{proposition} \label{prop_control_vertical_strip}
Assume that $q$ satisfies the assumptions of  \cref{th-main-pos-intro}. Let $\omega = (a,b)\times \T$, with
$-L_- \leq a < b \leq L_+$. Then 
\begin{itemize}
\item if $0< a$, the Baouendi-Grushin equation is null-controllable on $\o$ in time $T> \dagmon(a)/q'(0)$,
\item if $b< 0$, the Baouendi-Grushin equation is null-controllable on $\o$ in time $T> \dagmon(b)/q'(0)$,
\item if $a<0<b$, the Baouendi-Grushin equation is null-controllable on $\o$ in time $T>0$.
\end{itemize}
\end{proposition}

\begin{proof}[Proof of  \cref{th-main-pos-intro}]

\begin{figure}[h]
\centering
\input{omega0_left_right.tex}
\caption{Definition of $\omega_-$ (red) and $\omega_+$ (blue).}
\label{fig:omega 0 left right}
\end{figure}
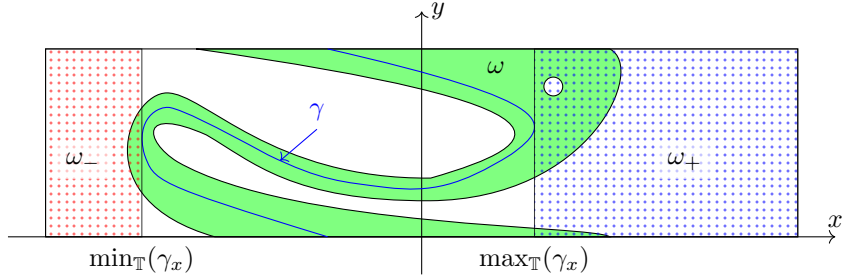
We set $\omega_- = (-L_-,\min_\T\left(\gamma_x\right)) \times \T$ and $\omega_+ = (\max_\T (\gamma_x),L_+) \times \T$ (see \cref{fig:omega 0 left right}). \Cref{prop_control_vertical_strip}
implies that the Baouendi-Grushin equation is null-controllable either on $\o_-$ or on $\o_+$ in any time $T$ such that 
\[
T> \frac 1 {q'(0)} \max\big(\dagmon(\min_\T\left(\gamma_x\right)), \dagmon(\max_\T (\gamma_x)) \big).
\]
Consequently, for any initial condition $f_0 \in L^2(I\times \T)$, there exist $u_-\in L^2((0,T);\, L^2(\o_-))$ and $u_+\in L^2((0,T);\, L^2(\o_+))$ such that
$f_-$ and $f_+$ solutions of
\[
\left\{\begin{array}{ll}
         (\partial_t-\partial_x^2 - q(x)^2 \partial_y^2)f_\pm(t,x,y) = \mathds 1_{\o_\pm} u_\pm(t,x,y),&\quad t\in (0,T), x\in I, y\in \T;\\
         f_\pm(t,x,y) = 0,&\quad t\in (0,T), x\in\partial I, y\in \T; \\
         f_\pm(0,x,y) = f_0(x,y) , & \quad x \in I, \ y \in \T;
        \end{array}\right.
\]
satisfy $f_\pm(T,\cdot,\cdot) = 0$ in $I\times \T$.

By definition of $\gamma$, $\omega_+$ and $\omega_-$ are included in two distinct connected components of $(I \times \T) \setminus \gamma(\T)$. 
As a consequence, we can construct $\chi \in C^\infty(I\times \T)$  such that
$\chi \equiv 1$ in $\omega_+$,
$\chi \equiv 0$ in $\omega_-$ and $\supp(\nabla \chi) \subset \o$ (see \cref{th-cutoff}).
Define $ f =  \chi f_- + (1-\chi) f_+$. It is easily verified that $f$ satisfies
\[
\left\{\begin{array}{ll}
         (\partial_t-\partial_x^2 - q(x)^2 \partial_y^2)f(t,x,y) = \mathds 1_{\o} u(t,x,y), &\quad t\in (0,T), x\in I, y\in \T;\\
         f(t,x,y) = 0,&\quad t\in (0,T), x\in\partial I, y\in \T; \\
         f(0,x,y) = f_0(x,y) , & \quad x \in I, \ y \in \T; \\
         f(T,x,y) = 0 , & \quad x \in I, \ y \in \T; 
        \end{array}\right.
\]
with a source term $u \in L^2((0,T);\, L^2(\o))$.
\end{proof}

\section{Lack of null-controllability}\label{sec_negative_result}
In this section, we prove the following case of \cref{th-main-neg-intro}.
\begin{theorem}\label{th-main-neg}
Assume that $q\in C^2(\overline I)$ is such that $q(0) = 0$, $q'(0) > 0$ and $q(x) \neq 0$ whenever $x\neq 0$.
Let $\o$ be an open subset of $I\times \T$. Assume that there exist $a \in (-L_-,0)$ and $y_0 \in \T$ such that 
\[
\distance\big((a,L_+)\times \{y_0\}, \o\big) > 0.
\]
Then, the generalized Baouendi-Grushin equation~\eqref{eq-grushin-gen-contr} is not null-controllable on $\o$ in time $T< \dagmon(a)/q'(0)$.
\end{theorem}

\begin{remark}\label{rk-main-neg}
By changing $x$ in $-x$, $I$ in $-I$ and $q$ in $-q$, this theorem implies that if $b\in(0,L_+)$ and if $\distance\big((-L_-,b)\times \{y_0\}, \o\big) > 0$, then the generalized Baouendi-Grushin equation~\eqref{eq-grushin-gen-contr} is not null-controllable on $\o$ in time $T< \dagmon(b)/q'(0)$.

To completly prove \cref{th-main-neg-intro}, there are two more cases: 
\begin{itemize}
\item $\distance\big((a,b)\times\{y_0\}, \o\big)>0$ with $-L_-< a <0<b<L_+$, lack of null-controllability in time $T< \min(\dagmon(a),\dagmon(b))/q'(0)$;
\item $\distance\big((-L_-,L_+)\times\{y_0\}, \o\big)>0$, lack of null-controllability in any time $T>0$.
\end{itemize}
The proofs of these cases are minor modifications of the one of \cref{th-main-neg}. We mention in footnotes the most important modifications and leave the details to the reader.
\end{remark}

\begin{figure}[h]
\centering
\input{omega_neg}
\caption{In green, the domain $\omega$. If a horizontal segment stays at positive distance from $\omega$, it can be thickened into a rectangle that is disjoint from $\omega$.}
\label{fig:omega_neg}
\end{figure}
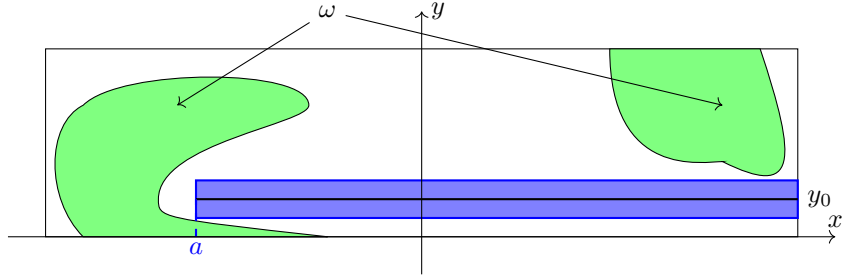

Under the hypotheses of this theorem, there exists a closed interval $W_0$ that is a neighborhood  of $y_0$ and such that $\o \cap \big([a,L_+)\times W_0\big) = \emptyset$ (see \cref{fig:omega_neg}). To prove \cref{th-main-neg}, we assume without loss of generality that $\o$ is the complement of the rectangle $[a,L_+)\times W_0$: 
\begin{equation}\label{eq-def-omega}
\o = (I\times \T)\setminus \big([a,L_+)\times W_0\big).
\end{equation}

\subsection{Observability inequality}
Using standard duality arguments (see~\cite[theorem~2.44]{Coron07}), the null-controllability of the generalized Baouendi-Grushin equation~\eqref{eq-grushin-gen-contr} is equivalent to the following \emph{observability inequality}: there exists $C>0$ such that for every $g_0\in L^2(I\times \T)$, the solution $g$ of
\begin{equation}
 \label{eq-grushin-gen}
 \left\{\begin{array}{ll}
         (\partial_t-\partial_x^2 - q(x)^2 \partial_y^2)g(t,x,y) = 0,&\quad t\in (0,T), x\in I, y\in \T;\\
         g(t,x,y) = 0,&\quad t\in (0,T), x\in \partial I, y\in \T;\\
         g(0,x,y) = g_0(x,y)&\quad (x,y) \in I\times \T,
        \end{array}\right.
\end{equation}
satisfies
\begin{equation}
 \label{eq-obs}
 \|g(T,\cdot,\cdot)\|_{L^2(I\times \T)}^2 \leq C \|g\|_{L^2((0,T)\times \o)}^2.
\end{equation}
To prove \cref{th-main-neg}, we proceed in two steps: we prove that the observability inequality~\eqref{eq-obs} implies an inequality on polynomials, and then we disprove this new inequality.\footnote{Actually, we could reformulate this proof to directly construct a counterexample to the observability inequality~\eqref{eq-obs}. 
}

\subsection{Model case}
We start with a model equation, that we study to showcase the main ideas of the proof of \cref{th-main-neg} without some of the more technical aspects. Consider the Baouendi-Grushin equation on $\R\times \T$:
\begin{equation}
 \label{eq-grushin-R}
 (\partial_t -\partial_x^2 - x^2\partial_y^2)g(t,x,y) = 0.
\end{equation}
	Let $\omega\subset \R \times \T$ be open. We say that the Baouendi-Grushin equation \eqref{eq-grushin-R} is observable on 
	$\o$ in time $T>0$ if there exists $C>0$ such that for all $g$
	solution of \eqref{eq-grushin-R}, the following observability
	inequality holds:
	\begin{equation} \label{eq-obs-R}
	\Vert g(T,\cdot,\cdot)\Vert^2_{L^2(\R\times \T)}
	\leq C \Vert g\Vert_{L^2((0,T)\times \o)}^2.
	\end{equation}
	
We prove the following theorem.
\begin{theorem}\label{th-model-neg}
let $a>0$, $W_0\subset \T$ a closed interval with non-empty interior and  
\begin{equation}\label{eq-def-o-model}
\o = (\R\times \T)\setminus([-a,a]\times W_0).
\end{equation}
Let $T>0$ such that
\[
T < \frac{a^2}2.
\]
The Baouendi-Grushin equation~\eqref{eq-grushin-R} is not observable on $\o$ in time $T$.
\end{theorem}

Before going into the proof, let us examine some solutions of the Baouendi-Grushin equation that are concentrated around $x=0$. Taking the $n$-th Fourier coefficient in $y$ of $g$, which we will denote by $\hat g(t,x,n)$, we get
\begin{equation}
 \label{eq-grushin-R-fourier}
 (\partial_t -\partial_x^2 +n^2x^2)\hat g(t,x,n) = 0.
\end{equation}
Thus, the Baouendi-Grushin equation is transformed into a family of parabolic equations $(\partial_t + H_n)g_n = 0$, where $H_n$ is the \emph{harmonic oscillator} $-\partial_x^2 + n^2 x^2$. The spectral properties of the harmonic oscillator are well-known (see, e.g., \cite[\S1.3]{Helffer13} or \cref{sec:Davies-operator}), and in particular the first eigenvalue is $|n|$ with associated eigenfunction $\f_n(x) = (|n|/\pi)^{1/4} \eu^{-|n|x^2\!/2}$. Thus, if $(a_n)_{n>0}$ is a complex-valued sequence with only a finite number of nonzero terms, the function $g$ defined by 
\begin{equation}\label{eq-model-solution}
 g(t,x,y) \coloneqq \sum_{n>0} a_n \eu^{\iu n y - n x^2\!/2 - nt}
\end{equation}
is a solution of the Baouendi-Grushin equation~\eqref{eq-grushin-R}. We will look for a counterexample of the observability inequality~\eqref{eq-obs-R} in this class of functions.

This solution can be written as $g(t,x,y) = g_\pol(\eu^{\iu y - t -x^2\!/2})$ with
\begin{equation}\label{eq-model-polyn}
 g_\pol(z) \coloneqq \sum_{n>0} a_n z^n.
\end{equation}
We will use the fact that $g$ is a polynomial in $z = \eu^{\iu y - t -x^2\!/2}$ to rewrite the observability inequality we want to disprove as an inequality on polynomials. More precisely, we have the following estimate.
\begin{lemma}\label{th-model-polyn-est}
Assume that the observability inequality for the Baouendi-Grushin equation \eqref{eq-obs-R} holds. Let $U\subset \C$ be defined by (see \cref{fig:U})
\[
U = D(0, \eu^{-a^2\!/2}) \cup 
    \{ z\in \C\colon  |z|<1,\ \arg(z)\notin W_0 \}. 
\]
Then, there exists $C>0$ such that for every polynomial $p\in \C[X]$,
\[
\|p\|_{L^2(D(0,\eu^{-T}))} \leq C \|p\|_{L^\infty(U)}.
\]
\end{lemma}
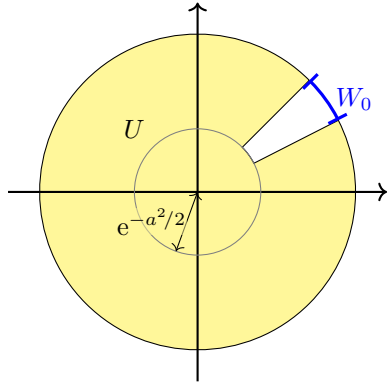
\begin{figure}
\centering
\input{U}
\caption{The domain $U$.}
\label{fig:U}
\end{figure}%
\begin{proof}
\step{Observability inequality}
Let $p(z) = \sum_{n\geq 0} a_n z^n$ a polynomial and set $g_\pol(z) = zp(z) = \sum_{n>0} a_{n-1}z^n$. The discussion above shows that $g$ defined by
\[
g(t,x,y) = g_\pol(\eu^{\iu y -t - x^2\!/2})
\]
is a solution of the Baounedi-Grushin equation \eqref{eq-grushin-R}. The observability inequality on this class of functions reads
\begin{equation}
\label{eq-obs-grushin-R}
\int_{\R\times \T} |g(T,x,y)|^2 \diff x\diff y \leq 
C\int_{[0,T]\times \o} |g_\pol(\eu^{\iu y -t - x^2\!/2})|^2 \diff t\diff x\diff y.
\end{equation}

\step{Left-hand side of the observability inequality \eqref{eq-obs-grushin-R}}
Since the functions $\psi_n\colon(x,y)\mapsto \eu^{\iu n y - n x^2\!/2}$ are orthogonal in $L^2(\R\times \T)$, the left-hand side can we rewritten as
\begin{align*}
\int_{\R\times \T} |g(T,x,y)|^2 \diff x\diff y 
&=\int_{\R\times \T} \Big|\sum_{n>0} a_{n-1} \eu^{-nT} \psi_n(x,y)\Big|^2 \diff x\diff y\\
&= \sum_{n>0} |a_{n-1}|^2 \eu^{-2nT} \|\psi_n\|_{L^2(\R\times\T)}^2\\
&= \sum_{n>0}  \frac{2\pi^{3/2}}{\sqrt{n}} |a_{n-1}|^2 \eu^{-2nT}.
\end{align*}
Elementary computations in polar coordinates prove that the functions $z\mapsto z^n$ are orthogonal in $L^2(D(0,R),m)$, where $m$ is the Lebesgue measure on $\C \simeq \R^2$,
and that for $R>0$
\[
\|z^n\|_{L^2(D(0,R),m)}^2 = \frac {\pi R^{2n+2}}{n+1}.
\]
Thus,
\begin{align}
\|p\|_{L^2(D(0,\eu^{-T}))}^2 &= \sum_{n\geq 0} \frac{\pi}{n+1} |a_n|^2 \eu^{-2(n+1)T}\notag\\
&\leq \sum_{n > 0} \frac\pi{\sqrt n} |a_{n-1}|^2 \eu^{-2nT}\notag\\
&= \frac1{2\sqrt \pi} \int_{\R\times \T} |g(T,x,y)|^2 \diff x\diff y.\label{eq-obs-R-lhs-bound}
\end{align}

\step{Right-hand side of the observability inequality \eqref{eq-obs-grushin-R}}
We write the right-hand side of the observability inequality by making the change of variables $(x,z) = (x,\eu^{-t+\iu y -x^2\!/2})$. We have $\diff x \diff m( z) = |z|^2 \diff t\diff x\diff y$. Thus, if we denote by $\Omega \subset \R\times \C$ the image of $(0,T)\times \o$ by this change of variables, we have
\begin{align}
\int_{(0,T)\times \o} |g_\pol(\eu^{-t + \iu y - x^2\!/2})|^2 \diff t \diff x\diff y
&= \int_{\Omega} |g_\pol(z)|^2 |z|^{-2} \diff x \diff m(z)\notag\\
&= \int_{\Omega} |p(z)|^2 \diff x \diff m(z).\label{eq-obs-R-rhs-est-1}
\end{align}
By definition, $\Omega = \bigcup_{x\in \R} (\{x\}\times \Dc_x)$, where
\begin{equation*}\label{eq-def-Dcx}
\Dc_x = \{\eu^{-t+\iu y -x^2\!/2},\ t\in (0,T), (x,y)\in \o\}.
\end{equation*}

We claim that for every $0<t<T$ and $x\in \R$, $\Dc_x \subset U$ (this is the reason we defined $U$ this way). Indeed, if $-a<x<a$, and $(x,y)\in \o$, then, by definition of $\o$ as the complement of $[-a,a]\times W_0$ (\cref{eq-def-o-model}), we necessarily have $y\notin W_0$. It follows that $\eu^{-t+\iu y -x^2\!/2}\in U$. In the case $x\notin [-a,a]$, we have $x^2\!/2 > a^2\!/2$. Then, $|\eu^{-t+\iu y -x^2\!/2}| < \eu^{-a^2\!/2}$. It follows again that $\eu^{-t+\iu y -x^2\!/2}\in U$. Thus, using Hölder inequality in \cref{eq-obs-R-rhs-est-1},
\begin{align}
\int_{(0,T)\times \o} |g_\pol(\eu^{-t + \iu y - x^2\!/2})|^2 \diff t\diff x\diff y
&= \int_{x\in\R} \int_{\Dc_x} |p(z)|^2 \diff m(z) \diff x\notag\\
&\leq \int_{x\in \R} m(\Dc_x) \|p\|_{L^\infty(U)}^2 \diff x.\notag
\intertext{Since $\Dc_x \subset D(0,\eu^{-x^2\!/2})$, $m(\Dc_x) \leq \pi \eu^{-x^2}$. Hence}
\int_{(0,T)\times \o} |g_\pol(\eu^{-t + \iu y - x^2\!/2})|^2 \diff t\diff x\diff y
&\leq (\pi)^{3/2} \|p\|^2_{L^\infty(U)}. \label{eq-obs-R-rhs-bound}
\end{align}

Now, plugging the lower-bound of the left-hand side~\eqref{eq-obs-R-lhs-bound} and the upper bound of the right-hand side~\eqref{eq-obs-R-rhs-bound} into the observability inequality~\eqref{eq-obs-grushin-R}, we obtain a constant $C>0$ such that
\[
\|p\|_{L^2(D(0,\eu^{-T}))}^2 \leq C \|p\|_{L^\infty(U)}^2. \qedhere
\]
\end{proof}

\begin{figure}
 \begin{minipage}[c]{0.5\textwidth}
\centering
\input{runge-R}
 \end{minipage}\hfill%
 \begin{minipage}[c]{0.5\textwidth}
\caption{When the disk $D(0,\eu^{-T})$ (in red) is not included in $U$, we can find holomorphic functions that are small in $U$ but arbitrarily large in $D(0,\eu^{-T})$. For instance, we can construct with Runge's theorem a sequence of polynomials that converges to $z\mapsto(z-z_0)^{-1}$ away from the blue line.}
\label{fig:runge-R}
\end{minipage}
\end{figure}
\begin{proof}[Proof of \cref{th-model-neg}]
To disprove the inequality, we only have to disprove the inequality on polynomials given by the previous \cref{th-model-polyn-est}. If $T< a^2/2$, the disk $D(0,\eu^{-T})$ is not included in $U$ (see \cref{fig:runge-R}). For instance, if $y_1\in \mathring W_0$ and $\e>0$ is small enough, $z_0 = \eu^{\iu y_1 -T-\e}$ is not in $U$. In fact, the half-line $z_0[1,+\infty)$ stays at positive distance from $U$ (see \cref{fig:runge-R}). Then, according to Runge's theorem~\cite[theorem~13.9]{Rudin86}, there exists a sequence of polynomials $(p_k)$ that converges uniformly on every compact of $\C \setminus z_0[1,+\infty)$ to $z\mapsto (z-z_0)^{-1}$.

Since $\overline U$ is a compact subset of $\C \setminus z_0[1,+\infty)$, the sequence $p_k$ stays uniformly bounded on $U$, i.e., $\sup_k \|p_k\|_{L^\infty(U)} < +\infty$. But $z_0 \in D(0,\eu^{-T})$, and therefore $\|(z-z_0)^{-1}\|_{L^2(D(0,\eu^{-T}))} = +\infty$. Thanks to Fatou's lemma, this proves that $\|p_k\|_{L^2(D(0,\eu^{-T}))} \to +\infty$ as  $k\to + \infty$.

We have proved that $(p_k)$ is a counterexample to the inequality of \cref{th-model-polyn-est}, which concludes the proof of \cref{th-model-neg}.
\end{proof}

\subsection{From the model case to the generalized Baouendi-Grushin equation}\label{sec-sketch-gen} Now, our goal is to adapt the strategy used in the model case to the generalized Baouendi-Grushin equation~\eqref{eq-grushin-gen-contr}.
In the generalized Baouendi-Grushin equation, if we take the $n$-th Fourier coefficient in $y$ of $g$, we get
\begin{equation}
 \label{eq-grushin-gen-fourier}
 (\partial_t -\partial_x^2 +n^2q(x)^2)\hat g(t,x,n) = 0.
\end{equation}
Recall that for $n\geq 0$, $P_n$ is the unbounded operator $-\partial_x^2 + n^2 q^2$ on $L^2(I)$ with Dirichlet boundary conditions. We will denote by $\l_n$ the first eigenvalue of $P_n$ and by $\f_n$ a corresponding eigenfunction. Notice that $\f_n$ is not required to be normalized in $L^2(I)$. Then, we will look for a counterexample of the observability inequality~\eqref{eq-obs} with solutions of the generalized Baouendi-Grushin equation~\eqref{eq-grushin-gen} of the form
\begin{equation}
 \label{eq-obs-counter}
 g(t,x,y) \coloneqq \sum_{n\geq 0} a_{n} \f_n(x) \eu^{\iu n y - \lambda_n t}.
\end{equation}
Heuristically, this should work because we expect the eigenfunction $\f_n$ to be localized around $x = 0$ as $n\to +\infty$, in which case the operator $-\partial_x^2 + n^2 q^2$ looks like $-\partial_x^2 +n^2 q'(0)^2 x^2$, and the eigenvalue and eigenfunction look like $\l_n \approx nq'(0)$ and $\f_n(x) \approx n^{1/4} \eu^{-nq'(0)x^2\!/2}$. So the solutions $g$ defined above look like the solutions used to treat the model case (\cref{eq-model-solution}), up to a factor $q'(0)$.

In fact, a better approximation of $\f_n$ would be the so-called WKB approximation\footnote{We write here the first term in the WKB expansion of eigenfunctions, and only in dimension $1$, because it is enough for our purposes. But such a construction can be refined with more terms and in higher dimension \cite[Chapter 3 \&\ Chapter 6, theorem A.3]{DS99}.

Also, in the differential equation that defines $c_0$, we divide by $\dagmon'(x)$, which is equal to $0$ at $x = 0$. But the numerator is also $0$ at $x = 0$, and simple Taylor expansions at $x=0$ proves that the quotient appearing in the differential equation for $c_0$ is actually well-defined at $x=0$.}
\begin{equation}\label{eq-WKB-selfadjoint}
\begin{array}{l}
\f_n(x) \approx n^{1/4} c_0(x) \eu^{-n \dagmon(x)},\\
c_0'(x) = \frac{q'(0) - \dagmon''(x)}{2\dagmon'(x)}c_0(x),\\
c_0(0) = 1.
\end{array}
\end{equation}

Thus, we have 
\begin{equation}
 \label{eq-obs-counter-approx}
 g(t,x,y) \approx  c_0(x)\sum_n a_{n} \eu^{n(\iu y - q'(0)t - \dagmon(x))},
\end{equation}
i.e., $g$ can almost be written as $g(t,x,y) \approx c_0(x)g_\pol(\eu^{\iu n y -t q'(0) -\dagmon(x)})$, where $g_\pol$ is the polynomial
\[
g_\pol(z) \coloneqq \sum_n a_{n} z^n.
\]

Let us write this in an exact way. Consider $\tilde \f_n(x) \coloneqq n^{1/4} \eu^{-nq'(0)x^2\!/2}$ and let $\Pi_n$ be the spectral projection associated to the first eigenvalue $\l_n$ of $P_n$. We define\footnote{We could also have chosen $\tilde \f_n$ to be the WKB expansion defined previously, which would be a better approximation of the eigenfunction. But since we are projecting on the actual eigenfunction afterwards, this is not necessary.}
\begin{equation}
\label{eq-def-vn}
\f_n \coloneqq \Pi_n \tilde \f_n.
\end{equation}
We will see later that $\f_n\neq 0$, at least if $n$ is large enough. Let $\e\in (0,1)$, that we need for technical reason, and that we will later choose close to 0. We define $\g_{t,x}(n)$ by
\begin{equation}\label{eq-def-gamma-tx}
\g_{t,x}(n-1) \coloneqq \eu^{-t(\l_n - q'(0)n)} \f_n(x)\eu^{n\dagmon(x)(1-\e)}.
\end{equation}
The shift of $n$ in the definition is linked to the fact that we will consider $p(z) = g_\pol(z)/z$, as we did in the model case. Then, the solution $g$ defined in~\cref{eq-obs-counter} can be written as
\begin{equation}
 \label{eq-obs-counter-exact}
 g(t,x,y) =  \sum_n a_{n} \g_{t,x}(n-1)\eu^{n\big(\iu y - q'(0)t - (1-\e)\dagmon(x)\big)}.
\end{equation}
In some sense, this formula tells us that $g$ can be written as ``pseudo-differential-type'' operator applied to the ``model solution'' $g_\pol(\eu^{\iu y - q'(0)t -(1-\e)\dagmon(x)})$. To successfully adapt the strategy used for the model Baouendi-Grushin equation, we need some continuity estimates for these ``pseudo-differential-type'' operators. We claim that the following estimate holds.
\begin{lemma}\label{th-pseudo-polyn-est}
Let $T>0$ and $\e>0$. Define $\g_{t,x}$ as in \cref{eq-def-gamma-tx}. Let $\Op{\g_{t,x}}$ be the operator on polynomials defined by
\[
\Op{\g_{t,x}}\Big(\sum a_n z^n\Big) = \sum \g_{t,x}(n) a_n z^n.
\]
Let \(X\noic\) be a compact subset of $\C$. Let \(V\noic\) be an open neighborhood of $X$ that is star-shaped with respect to $0$. There exist $C>0$ and $N\in\N$ such that for every polynomial $p \in \C[X]$ with a zero of order $N$ at $0$ and for every $0<t<T$ and $x\in I$,
\[
\lVert\Op{\g_{t,x}}(p)\|_{L^\infty(X)} \leq C\|p\|_{L^\infty(V)}.
\]
\end{lemma}

As $\g_{t,x}(n)$ is related to the eigenvalues and eigenfunctions of $P_n$, proving this lemma requires a spectral analysis of this operator. What is more surprising is that we actually need a spectral analysis of $P_\nu$ \emph{when $\nu$ is not necessarily real}, meaning we have to do some \emph{nonselfadjoint} spectral analysis. 
We will prove \cref{th-pseudo-polyn-est} in \cref{sec-pseudo-polyn-est} with the spectral analysis done in the rest of \cref{sec-spectral} and a general estimate on operators on polynomials~\cite[theorem~18]{Koenig17}.

We will also use the relatively elementary bounds on $\l_n$ and $\|\f_n\|_{L^2(I)}$ given by the following proposition:
\begin{proposition}\label{th-lower-eigen}
In the limit $n\to + \infty$, $\l_n = nq'(0) + o(n)$. Moreover, there exist $c>0$ and $N\geq 0$ such that for every $n\geq N$, $\|\f_n\|_{L^2(I)} \geq c$.
\end{proposition}
This proposition is standard (see, e.g., \cite[theorem 4.23 \&\ Eq.~(4.20)]{DS99}), nevertheless, for the reader convenience, we provide a proof in \cref{sec-agmon}.

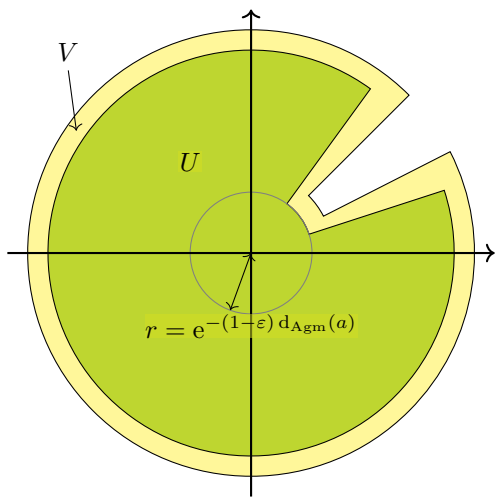
\begin{figure}
\centering
\input{D}
\caption{The domains $U$ and $V$.}
\label{fig:U-gen}
\end{figure}%
With these two estimates, we prove the following version of \cref{th-model-polyn-est} adapted for the generalized Baouendi-Grushin equation.
\begin{lemma}\label{th-polyn-est}
Assume that the observability inequality~\eqref{eq-obs} for the generalized Baouendi-Grushin equation holds. Let $\e>0$ and let $U\subset \C$ be defined by (see \cref{fig:U-gen})\footnote{In the variant of \cref{th-main-neg} where $\distance((a,b)\times\{y_0\},\o)>0$ mentioned in \cref{rk-main-neg}, we have to add $D(0,\eu^{-(1-\e)\dagmon(b)})$ to $U$. In the variant where $\distance(I\times\{y_0\},\o)>0$, $U$ is only the pacman $\{ z\in \C\colon  \arg(z)\notin W_0 \}$. Their proofs are minor adaptations and are left to the reader.\label{foot-variant1}}
\[
U = D(0, \eu^{-(1-\e)\dagmon(a)}) \cup 
    \{ z\in \C\colon |z|<1, \,  \arg(z)\notin W_0 \}. 
\]
Let \(V\noic\) be a neighborhood of \(\overline{U}\noic\) that is star-shaped with respect to $0$. Then, there exist $C>0$ and $N\in \N$ such that for every polynomial $p\in \C[X]$ with a zero of order $N$ at $0$, we have
\[
\|p\|_{L^2(D(0,\eu^{-q'(0) T(1+\e)}))} \leq C \|p\|_{L^\infty(V)}.
\]
\end{lemma}
\begin{proof}
The proof mostly follows the one of \cref{th-model-polyn-est}, but with the error term $\g_{t,x}$ which will be handled by \cref{th-pseudo-polyn-est}. Let $N>0$ as in \cref{th-lower-eigen} and \cref{th-polyn-est}. Let $p(z) = \sum_{n\geq N} a_n z^n$ a polynomial and $g_\pol(z) = zp(z)$. The discussion above shows that $g$ defined by
\begin{equation}\label{eq-def-g-gen}
\begin{aligned}
g(t,x,y) &= \sum_{n > N} a_{n-1} \f_n(x) \eu^{-\l_n t + \iu n y}\\ 
&=\sum_{n > N} a_{n-1} \g_{t,x}(n-1) \eu^{n\big(\iu y -q'(0)t - (1-\e)\dagmon(x)\big)}
\end{aligned}
\end{equation}
is a solution of the Baouendi-Grushin equation \eqref{eq-grushin-gen}.

\step{Left-hand side of the observability inequality \eqref{eq-obs}}
Since the functions $\psi_n\colon(x,y)\mapsto \f_n(x)\eu^{\iu n y}$ are orthogonal, the left-hand side can we rewritten as
\begin{align*}
\int_{I\times \T} |g(T,x,y)|^2 \diff x\diff y 
&=\int_{I\times \T} \Big|\sum_{n\geq N} a_{n-1} \eu^{-\l_n T} \psi_n(x,y)\Big|^2 \diff x\diff y\\
&= \sum_{n\geq N} |a_{n-1}|^2 \eu^{-2\l_n T} \|\psi_n\|_{L^2(I\times \T)}^2,
\intertext{using the lower bounds on $\|\f_n\|_{L^2(I)}$ given by \cref{th-lower-eigen}, we get that $\|\psi_n\|_{L^2(I\times \T)} \geq c > 0$ for $n\geq N$. Thus,}
\int_{I\times \T} |g(T,x,y)|^2 \diff x\diff y 
&\geq c \sum_{n\geq N} |a_{n-1}|^2 \eu^{-2\l_n T}.
\intertext{Now, thanks to the asymptotics for $\l_n$ given by \cref{th-lower-eigen}, there exists $C_\e> 0$ such that $\l_n \leq nq'(0)(1+\e) +C_\e$. Thus,}
\int_{I\times \T} |g(T,x,y)|^2 \diff x\diff y 
&\geq c\eu^{-2TC_\e} \sum_{n\geq N} |a_{n-1}|^2 \eu^{-2nq'(0)T(1+\e)}.
\end{align*}
As in the proof of \cref{th-model-polyn-est}, we denote by $m$ the Lebesgue measure on $\C \simeq \R^2$, the functions $z\mapsto z^n$ are orthogonal on $L^2(D(0,R),m)$ and $\|z^n\|_{L^2(D(0,R),m)}^2 = \pi R^{2n+2}/(n+1) $. Thus,
\begin{align}
\|p\|_{L^2(D(0,\eu^{-q'(0)T(1+\e)}))}^2 &= \sum_{n\geq  N} \frac{\pi}{n+1} |a_n|^2 \eu^{-2(n+1)q'(0)T(1+\e)}\notag\\
&\leq \sum_{n >  N} \pi |a_{n-1}|^2 \eu^{-2nq'(0)T(1+\e)}\notag\\
&\leq C \int_{I\times \T} |g(T,x,y)|^2 \diff x\diff y.\label{eq-obs-lhs-bound}
\end{align}

\step{Right-hand side of the observability inequality \eqref{eq-obs}}
We make the analogous change of variables as in the model case, but adapted to our case, i.e., $(x,z) = (x,\eu^{-q'(0)t+\iu y -(1-\e)\dagmon(x)})$. We have $\diff x \diff m(z) = q'(0)|z|^2 \diff t\diff x\diff y$. Thus, if we denote by $\Omega \subset I\times \C$ the image of $(0,T)\times \o$ by this change of variables, which is a subset of $I \times D(0,1)$, we have
\begin{align}
\int_{(0,T)\times \o} |g(t,x,y)|^2 \diff t \diff x\diff y
&= \frac1{q'(0)}\int_{\Omega} \bigg\lvert\sum_{n>N} a_{n-1} \g_{t,x}(n-1) z^n\bigg \rvert^2 |z|^{-2} \diff x \diff m(z)\notag\\
&=\frac1{q'(0)}\int_{\Omega} \left\lvert\g_{t,x}(z\partial_z)(p)(z)\right \rvert^2 \diff x \diff m(z).
\intertext{We keep for simplicity the notation $\g_{t,x}$ but of course $t$ is now a function of $(x,z)$. As in the model case, if $(x,z)\in \Omega$, then $z\in U$. Indeed, let $(x,z) \in \Omega$, i.e., $z = \eu^{-t+\iu y - \dagmon(x)(1-\e)}$ with $(x,y)\in \o$. If $a\leq x<L_+$, then, by definition of $\o$ as the complement of $[a,L_+)\times W_0$ (\cref{eq-def-omega}), we necessarily have $y\notin W_0$. It follows that $z\in U$. In the case $x< a <0$, since $\dagmon$ is decreasing on $[-L_-,0]$, we have $|z| = \eu^{-t-\dagmon(x)(1-\e)} < \eu^{-\dagmon(a)(1-\e)}$. It follows again that $z\in U$. Thus,}
\int_{(0,T)\times \o} |g(t,x,y)|^2 \diff t \diff x\diff y
&\leq C \sup_{(t,x)\in (0,T)\times I}\lVert\Op{\g_{t,x}}(p)\|^2_{L^\infty(U)}.\notag
\intertext{Now, we use the operator estimate of \cref{th-pseudo-polyn-est}, which gives}
\int_{(0,T)\times \o} |g(t,x,y)|^2 \diff t \diff x\diff y
&\leq C \| p\|^2_{L^\infty(V)}.\label{eq-obs-rhs-bound}
\end{align}

\step{Conclusion}
Now, plugging the lower-bound of the left-hand side~\eqref{eq-obs-lhs-bound} and the upper bound of the right-hand side~\eqref{eq-obs-rhs-bound} into the observability inequality~\eqref{eq-obs}, we get
\[
\|p\|_{L^2(D(0,\eu^{-q'(0)T(1+\e)}))}^2 \leq C \|p\|_{L^\infty(V)}^2. \qedhere
\]
\end{proof}

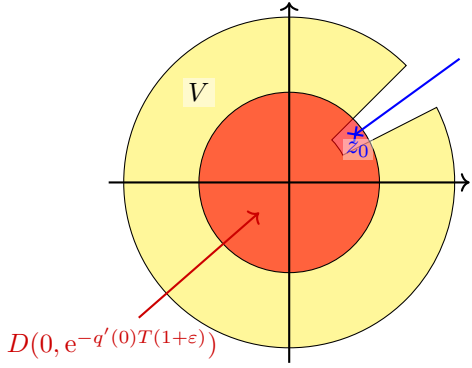
\begin{figure}
 \begin{minipage}[c]{0.5\textwidth}
\centering
\input{runge}
 \end{minipage}\hfill%
 \begin{minipage}[c]{0.5\textwidth}
\caption{When the disk $D(0,\eu^{-q'(0)T(1+\e)})$ (in red) is not included in $U$, we can find holomorphic functions that are small in $U$ but arbitrarily big in $D(0,\eu^{-q'(0)T(1+\e)})$. For instance, we can construct with Runge's theorem a sequence of polynomials that converges to $z\mapsto z^{N+1}(z-z_0)^{-1}$ away from the blue line.}
\label{fig:runge-gen}
\end{minipage}
\end{figure}

\begin{proof}[Proof of \cref{th-main-neg}]
As in the model case, we end the proof of non-null controllability by disproving the inequality on polynomials given by the previous \cref{th-polyn-est}. If 
\[
q'(0)T(1+\e)< \dagmon(a)(1-\e),
\]
the disk $D(0,\eu^{-q'(0)T(1+\e)})$ is not included\footnote{In the variant of \cref{th-main-neg} where $\distance((a,b)\times\{y_0\},\o)>0$ mentionned in \cref{rk-main-neg}, taking into account \cref{foot-variant1}, the condition becomes $q'(0)T(1+\e) < \min(\dagmon(a),\dagmon(b))(1-\e)$. In the variant where $\distance(I\times\{y_0\},\o)>0$, again taking into account \cref{foot-variant1}, $D(0,\eu^{-q'(0)T(1+\e)})$ is never included in $U$.} in $\overline U$, and we can chose a compact neighborhood $V$ of $\overline U$ that is star-shaped with respect to $0$ and such that $D(0,\eu^{-T(1+\e)})$ is not included in $\overline V$ (see \cref{fig:runge-gen}). Choose some $z_0\in D(0,\eu^{-q'(0)T(1+\e)})$ that is not in $\overline V$. Since $V$ is star-shaped with respect to $0$, the half-line $z_0[1,+\infty)$ stays at positive distance from $V$. Then, according to Runge's theorem~\cite[theorem~13.9]{Rudin86}, there exists a sequence of polynomials $(\tilde p_k)_k$ that converges uniformly on every compact of $\C \setminus z_0[1,+\infty)$ to $z\mapsto (z-z_0)^{-1}$. 

Set $p_k(z) \coloneqq z^{N+1} \tilde p_k$. We prove that $(p_k)$ is a counterexample to the inequality of \cref{th-polyn-est} with the same method as in the model case.
Since $\overline V$ is a compact subset of $\C \setminus z_0[1,+\infty)$, $p_k$ stays uniformly bounded on $V$, i.e., $\sup_k \|p_k\|_{L^\infty(V)} < +\infty$. But $z_0 \in D(0,\eu^{-q'(0)T(1+\e))})$, and therefore $\|z^{N+1}(z-z_0)^{-1}\|_{L^2(D(0,\eu^{-q'(0)T(1+\e)}))} = +\infty$. Thanks to Fatou's lemma, this proves that $\|p_k\|_{L^2(D(0,\eu^{-q'(0)T(1+\e)}))} \to +\infty$ as  $k\to + \infty$.

We have proved that the inequality of \cref{th-polyn-est} does not hold, which implies that the observability inequality~\eqref{eq-obs} does not hold either, which in turn implies that the generalized Baouendi-Grushin equation is not null-controllable. This holds for any $\e>0$ and any $T$ such that $(1+\e)q'(0)T< \dagmon(a)(1-\e)$. Thus, the generalized Baouendi-Grushin equation is not null controllable if $T< \dagmon(a)/q'(0)$.
\end{proof}

\section{Spectral Analysis}\label{sec-spectral}
As explained in \cref{sec-sketch-gen}, we need some spectral properties on the operator $P_\nu = -\partial_x^2 + \nu^2 q(x)^2$ with Dirichlet boundary conditions on $I$ (defined precisely in \cref{eq_def_Pnu}).
We start with an asymptotic of the first eigenvalue, and in following subsection, we prove some Agmon-type upper bound for the associated eigenfunctions.

For $\theta_0\in \left[0,\frac{\pi}{2} \right)$, we set
\begin{equation} \label{def:sector}
\sector \coloneqq \set{ \nu \in \mathbb{C} \colon  \abs \nu \geq 1 , \abs{\arg(\nu)} \leq \th_0} .
\end{equation}

\subsection{The first eigenvalue and corresponding spectral projection}

For $\b \in \C$ with $\Re(\b)>0$, we denote by $H_\b$ the non-selfadjoint harmonic oscillator $-\partial_x^2 + \b^2 x^2$ on $\R$. We refer to appendix \ref{sec:Davies-operator} for the precise definition and properties of $H_\b$.

In this paragraph we prove that the operator $P_\nu$ has an eigenvalue close to the eigenvalue $q'(0) \nu$ of the model operator $H_{q'(0) \nu}$, and that the corresponding spectral projection is also a perturbation of the spectral projection of $H_{q'(0)\nu}$. See proposition \ref{prop-lambda-n}.

For this we first prove that the resolvent of $P_\nu$ is a perturbation of the resolvent of $H_{q'(0)\nu}$, in the sense that the difference between these two resolvents is smaller than the resolvent of $H_{q'(0)\nu}$.

Notice that the resolvents of $P_\nu$ and $H_{q'(0)\nu}$ are not defined on the same space. We denote by $\1_I$ the operator that maps a function $v \in L^2(I)$ to its extension by 0 on $\R$.
Then $\1_I^*$ is the operator which maps $u \in L^2(\R)$ to its restriction on $I$: $\1_I^* u = u|_{I} \in L^2(I)$.

\begin{proposition}\label{prop:res-Pnu}
Let $\theta_0 \in \big[0,\frac \pi 2 \big)$. Let $\gamma > 0$ and $\e > 0$. For $\nu \in \sector$ we set (see \cref{def:Zc-beta})
\begin{equation} \label{def:Zc}
\tilde \Zc_\nu = \Zc_{q'(0)\nu,\e ,\gamma } =  \set{z \in \C \colon  \abs z \leq \gamma  q'(0) \abs \nu ,\ \distance(z , \Sp(H_{q'(0)\nu})) \geq \e q'(0) \abs \nu}.
\end{equation}
Then there exists $\nu_0 \geq 1$ such that for $\nu \in \sector$ with $\abs \nu \geq \nu_0$ and $z \in \tilde \Zc_\nu$ we have $z \in \rho(P_\nu)$, and 
\[
\sup_{z \in \tilde \Zc_\nu} \nr{(P_\nu-z)\inv - \1_I^* (H_{q'(0)\nu}-z)\inv \1_I }_{\Lc(L^2(I))} = \littleo[\nu \in \sector]{\abs \nu} {+\infty} \left( \frac 1 {\abs \nu} \right).
\]
\end{proposition}

\begin{proof} 
%\stepp 
For $\nu \in \Sigma_{\th_0}$ and $z \in \tilde \Zc_\nu$ we set 
\[
R_\nu(z) = \1_I^* (H_{q'(0)\nu}-z) \inv \1_I \quad \in \Lc(L^2(I)).
\]
Let 
\[
\rho \in \left]\frac 13,\frac 12 \right[.
\]
We consider a cut-off function $\chi \in C_0^\infty(\R,[0,1])$ supported in [-2,2] and equal to 1 on [-1,1]. Then for $\nu \in \Sigma_{\th_0}$ and $x \in \bar I$ we set 
\[
\chi_\nu(x) = \chi(\abs \nu^\rho x).
\]

\step{Approximation close to $x=0$} We first prove that if $\abs \nu$ is large enough then ${R_\nu(z) \chi_\nu (P_\nu-z)} - \chi_\nu$ extends to a bounded operator on $L^2(I)$ for all $z \in \tilde \Zc_\nu$, and 
\begin{equation} \label{estim-res-Rn}
\sup_{z \in \tilde \Zc_\nu} \nr{R_\nu(z) \chi_\nu (P_\nu-z) - \chi_\nu}_{\Lc(L^2(I))} \limt {\abs \nu} {+\infty} 0.
\end{equation}
Here and everywhere below it is implicitly understood that $\nu$ always belongs to $\sector$. 

Let $u \in \Dom(P_\nu)$. We have $\chi_\nu u \in \Dom(P_\nu)$ and, if $\abs \nu$ is large enough, $\1_I \chi_\nu u$ belongs to $\Dom(H_{q'(0)\nu})$. For $x \in \bar I$ we set 
\[
r(x) = q(x)^2 - q'(0)^2 x^2,
\]
so that, for $\abs \nu$ large enough,
\[
R_\nu(z) (P_\nu-z) \chi_\nu u = \chi_\nu u + \nu^2  R_\nu(z) r \chi_\nu u.
\]
The commutator $[\chi_\nu,P_\nu]$ of $\chi_\nu$ and $P_\nu$ is equal to $[\chi_\nu,P_\nu] = [\partial_x^2,\chi_\nu] = \chi_\nu'' + 2 \chi_\nu'\partial_x$, hence
\begin{align} 
R_\nu(z) \chi_\nu (P_\nu-z) u 
&= R_\nu(z) (P_\nu-z) \chi_\nu u + R_\nu(z) [\chi_\nu,P_\nu]u\notag\\
&= \chi_\nu u + \nu^2 R_\nu(z) r \chi_\nu u - R_\nu(z) \chi_\nu'' u + 2 R_\nu(z) (\chi_\nu' u)'.\label{eq-K-Rn}
\end{align}
By the resolvent estimate~\eqref{eq:res-Hb}, we have
\begin{equation} \label{estim:res-Hq}
\forall \nu \in \sector, \forall z \in \tilde \Zc_\nu, \quad \nr{(H_{q'(0) \nu} -z)\inv}_{\Lc(L^2(\R))} \leq \frac C {\abs \nu}.
\end{equation}
Since $\abs{r(x) \chi_\nu(x)} \lesssim \abs \nu^{-3\rho}$, this gives
\begin{equation*}
\nr{\nu^2 R_\nu(z) r \chi_\nu}_{\Lc(L^2(I))} \lesssim \abs \nu^{2 - 1 - 3\rho} \limt {\abs \nu} {+\infty} 0.
\end{equation*}
Similarly,
\[
\nr{R_\nu(z) \chi_\nu''}_{\Lc(L^2(I))} \lesssim \abs \nu^{2\rho -  1} \limt {\abs \nu} {+\infty} 0.
\]
Considering the last term in \cref{eq-K-Rn}, we have for $v \in L^2(\R)$ 
\begin{eqnarray} \label{eq:H-quad-form}
\lefteqn{\nr{ \partial_x \big(H_{q'(0)\bar \nu} - \bar z\big) \inv v }_{L^2(\R)}^2 + \bar \nu^2 q'(0)^2 \nr{x \big(H_{q'(0)\bar \nu} - \bar z\big)\inv v}_{L^2(\R)}^2}\\
\nonumber 
&& = \innp{H_{q'(0)\bar \nu} \big(H_{q'(0)\bar \nu} - \bar z\big)\inv v}{\big(H_{q'(0)\bar \nu} - \bar z\big)\inv v}_{L^2(\R)}\\
\nonumber
&& =  \innp{v}{\big(H_{q'(0)\bar \nu} - \bar z\big)\inv v}_{L^2(\R)}  + \bar z \nr{\big(H_{q'(0)\bar \nu} - \bar z\big)\inv v}_{L^2(\R)}^2.
\end{eqnarray}
We multiply by $\eu^{\iu\th}$ and take the real part. This gives, uniformly in $\nu \in \Sigma_{\th_0}$ and $z \in \tilde \Zc_\nu$,
\begin{equation} \label{eq:estim-partial-y-Hinv}
\nr{\partial_x \big(H_{q'(0)\bar \nu} - \bar z\big)\inv }_{\Lc(L^2(\R))}^2  \lesssim \frac {1}{\abs \nu} \quad \text{and} \quad \nr{x \big(H_{q'(0)\bar \nu} - \bar z\big)\inv}_{\Lc(L^2(\R))}^2 \lesssim \frac 1 {\abs \nu^3}.
\end{equation}
Taking the adjoint in the first inequality gives, for $\abs \nu$ large enough,
\begin{align*}
\nr{R_\nu(z) \partial_x (\chi_\nu' u)}_{L^2(I)} 
&\leq \nr{\big(H_{q'(0)\nu} -  z\big)\inv \partial_x (\1_I \chi_\nu' u)}_{L^2(\R)}\\
&\lesssim \abs \nu^{-\frac 12} \nr{\1_I \chi_\nu'u}_{L^2(\R)}\\
&\lesssim \abs\nu^{\rho-\frac 12} \nr{u}_{L^2(I)},
\end{align*}
and \cref{estim-res-Rn} follows.

\step{Approximation away from $x=0$}
There exists $c_0 \in (0,1]$ such that for all $x \in I$ we have 
\[
\abs{q(x)} \geq c_0 \abs x.
\]
On $L^2(I)$ we consider the operator 
\[
\tilde P_\nu = P_\nu + \nu^2 \1_{[-\abs \nu^{-\rho},\abs \nu^{-\rho}]},
\]
with domain $\Dom(\tilde P_\nu) = H^2(I) \cap H^1_0(I)$. It has compact resolvent, so its spectrum consists of eigenvalues. We have $q(x)^2 + \1_{[-\abs \nu^{-\rho},\abs \nu^{-\rho}]} \geq c_0^2 \abs \nu^{-2\rho}$ on $I$. Then, for $u \in \Dom(\tilde P_\nu)$ and $z \in \tilde \Zc_\nu$,
\[
\Re \big( \eu^{-\iu\th} \langle (\tilde P_\nu-z) u,u\rangle \big) \geq \big( \cos(\th) c_0^2 \abs \nu^{2 - 2\rho} - \gamma q'(0) \abs \nu \big) \nr{u}_{L^2(I)}^2.
\]
Thus, when $\abs \nu$ is so large that $\gamma q'(0) \abs \nu \leq \cos(\th_0) c_0^2 \abs \nu^{2 - 2\rho} / 2$ we have $\tilde \Zc_\nu \subset \rho(\tilde P_\nu)$ and, for $z \in \tilde \Zc_\nu$, 
\begin{equation} \label{eq:res-tilde-Pnu}
\big\| (\tilde P_\nu-z)\inv \big\|_{\Lc(L^2(I))} \leq \frac {2 \abs \nu^{2\rho-2}}{\cos(\th) c_0^2}.
\end{equation}
Then we have 
\begin{align*}
(\tilde P_\nu-z)\inv  (1-\chi_\nu) (P_\nu-z) u  
& = (\tilde P_\nu-z)\inv  (1-\chi_\nu) (\tilde P_\nu-z) u\\
& = (\tilde P_\nu-z)\inv  (\tilde P_\nu-z) (1-\chi_\nu) u + (\tilde P_\nu-z)\inv [1-\chi_\nu,\tilde P_\nu] u\\
& = (1-\chi_\nu) u - 2 (\tilde P_\nu-z)\inv (\chi_\nu' u)' + (\tilde P_\nu-z)\inv \chi_\nu'' u.
\end{align*}
As above we estimate 
\[
\big\| (\tilde P_\nu-z)\inv \chi_\nu'' \big\| _{\Lc(L^2(I))} \lesssim \abs \nu^{4\rho - 2} \limt {\abs \nu}{+\infty} 0
\]
and
\[
\big\| (\tilde P_\nu-z)\inv (\chi_\nu ' u)' \big\|_{L^2(I)} \lesssim \abs \nu^{\rho-1} \nr{\chi_\nu' u}_{L^2(I)} \lesssim \abs \nu^{2\rho -1} \nr{u}_{L^2(I)}.
\]
This proves that 
\begin{equation} \label{eq:estim-tilde-Pnu}
\sup_{z \in \tilde \Zc_\nu} \nr{(\tilde P_\nu-z)\inv  (1-\chi_\nu) (P_\nu-z) - (1-\chi_\nu)}_{\Lc(L^2(I))} \limt {\abs \nu}{+\infty} 0.
\end{equation}

\step{Conclusion} For $\nu \in \sector$ and $z \in \tilde \Zc_\nu$ we set 
\[
Q_\nu(z) = R_\nu(z) \chi_\nu + (\tilde P_\nu-z) \inv (1-\chi_\nu).
\]
By \cref{estim-res-Rn} and \cref{eq:estim-tilde-Pnu} we have for $u \in \Dom(P_\nu)$
\begin{eqnarray*}
\lefteqn{\nr{Q_\nu(z) (P_\nu-z)u - u}_{L^2(I)}} \\
&&\leq \nr{R_\nu(z) \chi_\nu (P_\nu-z) u - \chi_\nu u}_{L^2(I)} + \|(\tilde P_\nu-z)\inv (1-\chi_\nu) (P_\nu-z) u - (1-\chi_\nu) u \|_{L^2(I)} \\
&& = \littleo {\abs \nu}{+\infty} (1) \nr{u}_{L^2(I)}. 
\end{eqnarray*}
This proves that for $\abs \nu$ large enough, the operator $Q_\nu(z) (P_\nu-z) = 1 + (Q_\nu(z) (P_\nu-z) - 1)$ extends to a bounded operator on $L^2(I)$, which is invertible with inverse bounded uniformly in $\nu$, and 
\begin{equation} \label{eq:QP}
\nr{\big(Q_\nu(z) (P_\nu-z) \big)\inv - 1}_{\Lc(L^2(I))}\limt {\abs \nu}{+\infty} 0.
\end{equation}
In particular, $(P_\nu - z)$ is injective. Since it has compact resolvent, it is boundedly invertible, and
\[
(P_\nu-z)\inv = \big( Q_\nu(z) (P_\nu-z) \big)\inv Q_\nu(z).
\]
We get
\begin{equation} \label{eq:Pnu-Rnu}
\begin{aligned}
(P_\nu-z) \inv - R_\nu(z) 
& = \big( Q_\nu(z) (P_\nu-z) \big)\inv (\tilde P_\nu-z) \inv (1-\chi_\nu)\\
& - \big( Q_\nu(z) (P_\nu-z) \big)\inv R_\nu(z) (1-\chi_\nu)\\
& +  \Big(\big( Q_\nu(z) (P_\nu-z) \big)\inv - 1 \Big) R_\nu(z). 
\end{aligned}
\end{equation}
We prove that each term of the right-hand side is of size $O(\abs \nu^{-1})$.
For the first term we use \cref{eq:res-tilde-Pnu}. For the third we use \cref{eq:QP} and \cref{estim:res-Hq}. Finally, for the second term we observe that on $\supp(1-\chi_\nu)$ we have $\abs x \gtrsim \abs \nu^{-\rho}$ so for $u \in L^2(I)$ we have by the second inequality of \cref{eq:estim-partial-y-Hinv}
\begin{align*}
\big\| (1-\chi_\nu) \1_I^* (H_{q'(0)\bar \nu}-\bar z)\inv \1_I u\big\|^2_{L^2(I)}
 \lesssim \abs \nu^{2\rho} \big\|x (H_{q'(0)\bar \nu}-\bar z)\inv \1_Iu\big\|^2_{L^2(\R)} 
 \lesssim \abs \nu^{2\rho-3}\nr u^2_{L^2(I)}.
\end{align*}
Taking the adjoint gives
\[
\big\|\1_I^* (H_{q'(0)\nu} - z)\inv \1_I (1-\chi_\nu)\big\|_{\Lc(L^2(I))} = \littleo \nu {+\infty} \left( \frac 1 {\abs \nu} \right),
\]
and the conclusion follows from \cref{eq:Pnu-Rnu}.
\end{proof}

\begin{proposition} \label{prop-lambda-n}
Let $\th_0 \in \big[0,\frac \pi 2 \big)$. There exists $\nu_{\theta_0} \geq 0$ such that for $\nu \in \sector$ with $\abs \nu \geq \nu_{\theta_0}$ the operator $P_\nu$ has a unique eigenvalue $\lambda_\nu$ which satisfies $\abs{\lambda_\nu - q'(0)\nu} \leq q'(0) \abs \nu$. Moreover, $\lambda_\nu$ is algebraically simple and
\[
\abs{\l_\nu - q'(0) \nu} = \littleo[\nu \in \sector]{\abs \nu} {+\infty} \left(\abs \nu \right). 
\]
If we denote by $\Pi_\nu \in \Lc(L^2(I))$ the associated spectral projection, and by $\Pi^\h_{q'(0)\nu} \in \Lc(L^2(\R))$ the spectral projection associated to the eigenvalue $q'(0)\nu$ of $H_{q'(0)\nu}$, then
\[
\|\Pi_\nu - \1_I^* \Pi^\h_{q'(0)\nu}\1_I\|_{\Lc(L^2(I))} \limt[\nu \in \sector]{|\nu|}{+\infty} 0.
\]
\end{proposition}

\begin{proof}
 We recall that for $\nu \in \sector$ we have
\[
\Pi^\h_{q'(0) \nu} = \frac 1 {2\iu\pi} \int_{\abs{z - q'(0) \nu} = q'(0)|\nu|} (H_{q'(0)\nu}-z)\inv \diff z.
\]
Let $\nu_0$ be given by Proposition \ref{prop:res-Pnu} for $\e = 1$ and $\g = \frac 52$. If $\abs \nu \geq \nu_0$ we can set 
\[
B_\nu = \frac 1 {2\iu\pi} \int_{\abs{z - q'(0) \nu} = q'(0)|\nu|} (P_\nu-z)\inv \diff z.
\]
This is a projection of $L^2(I)$ whose range is the sum of the generalized eigenspaces of $P_\nu$ corresponding to the eigenvalues in the disk $D(q'(0) \nu, q'(0)\abs \nu)$. In particular the dimension of $\Ran(B_\nu)$ does not depend on $\nu \in \Sigma_{\theta_0}$ and is finite. We denote by $m$ this dimension and prove that $m=1$ by computing the trace $\Tr(B_\nu)$ of $B_\nu$.

By Proposition \ref{prop:res-Pnu} we have 
\[
\| B_\nu - \1_I^* \Pi^\h_{q'(0) \nu} \1_I \|_{\Lc(L^2(I))} \limt[\nu \in \sector]{|\nu|}{+\infty}  0.
\]

Let $F_\nu$ be a subspace of $L^2(I)$ of dimension $m+1$ which contains $\Ran(B_\nu)$ and $\1_I^* \f_{q'(0)\nu}^\h$. We consider an orthonormal basis $(e_{0,\nu},\dots,e_{m,\nu})$ of $F_\nu$. Then
\begin{align*}
\abs{\Tr\big(B_\nu - \1_I^* \Pi^\h_{q'(0) \nu} \1_I \big)} 
&\leq \sum_{j=0}^m \abs{\innp{\big(B_\nu - \1_I^* \Pi^\h_{q'(0) \nu} \1_I \big)e_{j,\nu}}{e_{j,\nu}}}\\
&\leq (m+1) \nr{B_\nu - \1_I^* \Pi^\h_{q'(0) \nu} \1_I}_{\Lc(L^2(I))} 
\limt[\nu \in \sector]{|\nu|}{+\infty} 0.
\end{align*}
Moreover, according to \cref{th-trace-eigenproj}, we have
\begin{equation*}
\Tr(\1_I^* \Pi^\h_{q'(0) \nu} \1_I)  \limt[\nu\in\sector]{|\nu|}{\infty} 1.
\end{equation*}
Thus,
\[
\Rank(B_\nu) = \Tr(B_\nu) \limt[\nu\in\sector]{|\nu|}{\infty} 1,
\]
and hence $m = 1$. This means that for $\nu \in \Sigma_{\theta_0}$ large enough the operator $P_\nu$ has a unique eigenvalue $\l_\nu$ in the disk $D(q'(0)\nu,q'(0)\abs \nu)$, and this eigenvalue is algebraically simple.

It remains to prove the estimate on $\l_\nu$. For this we reproduce the same argument with any $\e \in ]0,1]$. Then, given $\e \in ]0,1]$, there exists $\nu_\e \geq \nu_{\theta_0}$ such that the projection
\[
B_{\nu,\e} = \frac 1 {2\iu \pi} \int_{\abs{z - q'(0) \nu} = \e q'(0) |\nu|} (P_\nu-z)\inv \diff z
\]
is well defined and has rank 1. This implies that $\l_\nu$ belongs to $D(q'(0)\nu, \e q'(0) \abs\nu)$ when $\abs \nu \geq \nu_\e$ and concludes the proof.
\end{proof}

\subsection{Agmon Estimates of Eigenfunctions}\label{sec-agmon}
\begin{proposition}\label{th:agmon-eq}
Let $\nu \in \C$ such that $\Re(\nu)>0$. Let $\l$ be an eigenvalue of $P_\nu$ with associated eigenfunction $\f$. Set
\[
\mu\coloneqq \frac{\Re(\eu^{-\iu\arg(\nu)}\lambda)}{\cos(\arg(\nu))}.
\]
Let $\kappa \in C^2(I;\R) $  and $w \coloneqq \eu^{\nu \kappa}\f$. The following Agmon equality holds:
\[
\|w'\|_{L^2(I)}^2 + \int_I \left(|\nu|^2(q(x)^2 - \kappa'^2(x)) - \mu\right)|w(x)|^2 \diff x = 0.
\]
\end{proposition}

\begin{proof}
Let $\nu = |\nu|\eu^{\iu\th}$. The proof follows the usual Agmon's equality strategy. We have $\f = \eu^{-\nu \kappa} w$. Thus, using Leibniz' formula,
\begin{align*}
(P_\nu - \l) \f
 &=-(\eu^{-\nu \kappa})'' w - 2(\eu^{-\nu \kappa})' w' - \eu^{-\nu \kappa} w'' +(\nu^2 q^2-\l) \eu^{-\nu \kappa} w\\
 &=\left(-w'' +2\nu \kappa' w' +(-\nu^2 \kappa'^2 +\nu \kappa'' +\nu^2 q^2 - \l)w \right)\eu^{-\nu \kappa}.
\end{align*}
Since $(P_\nu - \l_\nu)\f = 0$, 
\begin{equation}
 -w'' +2\nu \kappa' w' +(-\nu^2 \kappa'^2 +\nu \kappa'' +\nu^2 q^2 - \l)w = 0.
\end{equation}
Multiplying by $\eu^{-\iu\th}\overline{w}$, integrating and taking the real part, we get
\begin{IEEEeqnarray*}{rCl} 
 0&=& \underbrace{-\Re\int_{\mathrlap I}\, \eu^{-\iu\th}w'' \overline{w}}_{I_1}
  +\underbrace{|\nu|\Re\int_{\mathrlap I}\, (2\kappa'w' \overline{w} + \kappa''|w|^2)}_{I_2} +\underbrace{\int_{\mathrlap I}\,\left( \cos(\th)|\nu|^2(q^2-\kappa'^2) - \Re(\eu^{-\iu\th}\l)\right) |w|^2}_{I_3}.
\end{IEEEeqnarray*}
Integrating by parts in $I_1$, we have $I_1 = \cos(\th) \|w'\|_{L^2(I)}^2$. 
Considering $I_2$, we get
\begin{equation*}
  I_2 =|\nu|\int_I(\kappa'|w|^2)' =0.
\end{equation*}
Thus,
\[
 0 = I_1 + I_3 =\cos(\th)\|w'\|^2_{L^2(I)} + \int_I\left( \cos(\th)|\nu|^2(q^2-\kappa'^2) - \Re(\eu^{-\iu\th}\l)\right) |w|^2,
\]
which is the claimed estimate multiplied by $\cos(\th)$.
\end{proof}

We will use \cref{th:agmon-eq} with $\kappa = (1-\e)\dagmon$, where $\dagmon$ defined in \cref{eq-def-agmon}. Up to this point, we assumed  $\f$ to be an eigenfunction of $P_\nu$, but we did not specified which one, neither how it is normalized. We do this in the following definition, which is the natural extension of the definition of $\f_n$ when $n\in\N$ (\cref{eq-def-vn}). For $\Re(\nu)>0$ that satisfies the hypotheses of~\cref{prop-lambda-n}, let $\tilde \f_\nu\ \in L^2(I)$ be defined by
\begin{equation*}
    \tilde \f_\nu(x) \coloneqq \nu^{1/4} \eu^{-q'(0) \nu  x^2\!/2},
\end{equation*}
and
\begin{equation}\label{eq-def-f-nu}
\f_\nu \coloneqq \Pi_\nu(\tilde \f_\nu),
\end{equation}
where $\Pi_\nu$ is the spectral projection for $P_\nu$ associated with $\lambda_\nu$, as defined in~\cref{prop-lambda-n}.

\begin{proposition}\label{th:est-vnu-l2}
Let $\f_\nu$ as in \cref{eq-def-f-nu}. Let $\theta_0 \in \big[0,\frac \pi 2 \big)$. There exists $C >0$ such that for every $\nu \in \sector$ with $|\nu|>1$,
\begin{gather*}
\|\tilde \f_\nu\|_{L^2(I)} \leq C;\\
\|\f_\nu\|_{L^2(I)} \leq C.
\end{gather*}
\end{proposition}

\begin{proof}
\step{Estimation on $\tilde \f_\nu$} Since $\tilde \f_\nu$ is a restriction of $\nu^{1/4}\eu^{-\nu q'(0) x^2 \!/2}$ on $\R$, we have
\begin{equation}\label{eq-est-vnu-2}
\|\tilde \f_\nu\|_{L^2(I)}^2 \leq |\nu|^{1/2}\int_\R \eu^{-2\Re(\nu) q'(0)x^2\!/2} \diff x
= \sqrt{\frac{\pi|\nu|}{q'(0)\Re(\nu)}} \leq \sqrt{\frac{\pi}{q'(0)\cos(\th_0)}}.
\end{equation}

\step{Estimation on $\f_\nu$} Using the notations of \cref{prop-lambda-n}, we have
\begin{equation}
\label{eq-est-vnu-1}
\|\f_\nu\|_{L^2(I)} \leq \Big(\|\1_I^* \Pi_{q'(0)\nu}^\h\1_I\|_{\Lc(L^2(I))} + \|\Pi_\nu - \1_I^*\Pi_{q'(0)\nu}^\h \1_I\|_{\Lc(L^2(I))}\Big)\|\tilde \f_\nu\|_{L^2(I)}.
\end{equation} Moreover, according to \cref{prop-lambda-n}, for $|\nu|$ large enough in $\sector$, 
\[
\|\Pi_\nu - \1_I^*\Pi_{q'(0)\nu}^\h\1_I\|_{\Lc(L^2(I))} \leq 1.
\]
Since the left-hand side is continuous in $\nu$, we have that for $|\nu|>1$ and $\nu\in \sector$,
\[
\|\Pi_\nu - \1_I^*\Pi_{q'(0)\nu}^\h\1_I\|_{\Lc(L^2(I))} \leq C.
\]
Finally, according to \cref{prop:Davies}, $\| \Pi^\h_{q'(0)\nu}\|_{\Lc(L^2)}$ stays bounded for $\nu\in\sector$, thus
\[
\|\1_I^*\Pi_{q'(0)\nu}^\h\1_I\|_{\Lc(L^2(I))} + \|\Pi_\nu - \1_I^*\Pi_{q'(0)\nu}^\h\1_I\|_{\Lc(L^2(I))} \leq C.
\]
Plugging this fact into~\cref{eq-est-vnu-1} and combined with the fact that $\|\tilde \f_\nu\|_{L^2(I)}$ is bounded (\cref{eq-est-vnu-2}), we get the claimed estimate.
\end{proof}

\begin{corollary}\label{prop-agmon-linfty}
Let $\f_\nu$ as in \cref{eq-def-f-nu}. Let $\theta_0 \in \big[0,\frac \pi 2 \big)$ and $\e \in (0,1)$. There exists $C>0$ such that for every $\nu \in \sector$,
\begin{gather*}
\int_I \big|\f_\nu(x) \eu^{\nu(1-\e)\dagmon(x)}\big|^2\diff x \leq C|\nu|.\\
\|\f_\nu \eu^{\nu(1-\e)\dagmon(x)}\|_{L^\infty(I)} \leq C|\nu|.
\end{gather*}

\end{corollary}
\begin{proof}
Let $w_\nu(x) \coloneqq \f_\nu(x) \eu^{\nu(1-\e)\dagmon(x)}$. According to the Agmon equality of \cref{th:agmon-eq} with $\f = \f_\nu$ and $\l = \l_\nu$, and denoting the corresponding $\mu$ by $\mu_\nu$, we have
\begin{equation}\label{eq:agmon-eq-l2}
\|w_\nu'\|_{L^2(I)}^2 + \int_I \big(|\nu|^2 q(x)^2(1-(1-\e)^2) - \mu_\nu\big)|w_\nu(x)|^2 \diff x = 0. 
\end{equation}
\step{First inequality}
Let $K>0$. We claim that if $x\in I$, $\nu \in \sector$ and $|\nu|\dagmon(x) > K$, then $|\nu|^2q(x)^2 > cK|\nu|$ for some $c$ depending on $q$, but not on $x\in I$, $\nu\in \sector$ nor $K$. 

Indeed, routine computations show that in the limit $x\to 0$, $q(x)^2 \sim q'(0)^2 x^2$ and $\dagmon(x) \sim q'(0)x^2\!/2$. Thus, $q(x)^2/\dagmon(x) \to 2q'(0)$ as $x \to 0$. Moreover, we assumed that for every $x\neq 0$, $q(x) \neq 0$, thus for every $x\neq 0$, $q(x)^2/\dagmon(x) \neq 0$. This proves the claim.

Thus, if $K$ is large enough, then for every $x\in I$ such that $\dagmon(x)|\nu| > K$
\[
|\nu|^2q(x)^2(1-(1-\e)^2) - \mu_\nu \geq c|\nu|,
\]
where we used the fact that $\mu_\nu \sim q'(0)|\nu|/\cos(\arg(\nu))$ because $\lambda_\nu = q'(0)\nu + o(\nu)$. Then, splitting the integral in the Agmon equality~\eqref{eq:agmon-eq-l2} into a part for $|\nu|\dagmon(x) < K$ and $|\nu|\dagmon(x) >K$, we get
\begin{align*}
c|\nu| \|w_\nu\|_{L^2(|\nu|\dagmon(x)>K)}^2
&\leq \int_{|\nu|\dagmon(x) < K} \big(\mu_\nu - |\nu|^2 q(x)^2(1-(1-\e)^2)\big)|w_\nu(x)|^2\diff x\\
&\leq C|\nu|^2 \|w_\nu\|_{L^2(|\nu|\dagmon(x) < K)}^2.
\intertext{We rewrite this as}
\|w_\nu\|_{L^2(|\nu|\dagmon(x)>K)}^2
&\leq C|\nu| \|w_\nu\|_{L^2(|\nu|\dagmon(x) < K)}^2.
\intertext{Adding $\|w_\nu\|_{L^2(|\nu|\dagmon(x) < K)}^2$ on each side, this proves that}
\|w_\nu\|_{L^2(I)}^2
&\leq C|\nu| \|w_\nu\|_{L^2(|\nu|\dagmon(x) < K)}^2.\\
\intertext{Using the definition of $w_\nu$, we see that for $|\nu|\dagmon(x) < K$, $|w_\nu(x)|\leq \eu^K |\f_\nu(x)|$. Thus, using also the property $\|\f_\nu\|_{L^2(I)}\leq C$ (\cref{th:est-vnu-l2})}
 \|w_\nu\|_{L^2(I)}^2
&\leq C|\nu| \|\f_\nu\|_{L^2(|\nu|\dagmon(x) < K)}^2 \leq C|\nu|.
\end{align*}

\step{Second inequality} We again use Agmon's equality~\eqref{eq:agmon-eq-l2} to get
\[
\|w_\nu'\|_{L^2(I)}^2 \leq \mu_\nu \|w_\nu\|_{L^2(I)}^2 \leq C|\nu|\|w_\nu\|_{L^2(I)}^2 \leq C|\nu|^2.
\]
The claimed estimate then follows from Sobolev's embedding of $H^1(I)$ into $L^\infty(I)$.
\end{proof}

We also prove the lower bound of \cref{th-lower-eigen} for $\f_n$ when $n\geq 0$ is large enough.\footnote{This theorem actually holds if $\nu$ ranges over $\sector$ by using the expression of $\Pi_\b^\h$. We don't need this, so we refrain from doing so.}
\begin{proof}[Proof of \cref{th-lower-eigen}]
The part about $\l_n$ was already proved in \cref{prop-lambda-n}. By definition of $\f_n$, we have
\[
\f_n = \Pi_n \tilde \f_n = \1_I^*\Pi^\h_{q'(0)n} \1_I \tilde \f_n + (\Pi_n - \1_I^*\Pi^\h_{q'(0)n}\1_I)\tilde \f_n.
\]
According to \cref{prop-lambda-n}, we have 
\[
\Pi_n - \1_I^*\Pi^\h_{q'(0)n}\1_I \limt{n}{+\infty} 0.
\]
Moreover, denoting by $\f_{\b,1}^\h(x) = (\Re(\b)/\pi)^{1/4}\eu^{-\b x^2\!/2}$ the first eigenvector of $H_\b$, we have for $\b>0$, $\Pi_\b^\h = \langle \f_{\b,1}^\h,\cdot\rangle \f_{\b,1}^\h$. Thus,
\begin{align*}
\langle \f_{q'(0)n,1}^\h, \1^*_I \tilde \f_n\rangle
&= \sqrt n \left(\frac{q'(0)}{\pi}\right)^{1/4} \int_I \eu^{- nq'(0)x^2} \diff x.
\intertext{The integral above is on $I$, but if we integrate on $\R$ instead, we only add a small error term. Thus,}
\langle \f_{q'(0)n,1}^\h, \1^*_I \tilde \f_n\rangle
&= \left(\frac\pi{q'(0)}\right)^{1/4} + \littleo{n}{+\infty}(1).
\end{align*}
Hence,
\[
\f_n = \left(\frac{\pi}{q'(0)}\right)^{1/4}\1_I^*\f_{q'(0)n,1}^\h + \littleo{n}{+\infty}(\|\tilde \f_n\|_{L^2(I)} +1).
\]
Since $\|\tilde \f_n\|_{L^2(I)}$ is bounded (\cref{th:est-vnu-l2}) and since $\|\1_I^* \f_{q'(0)n,1}^\h\|_{L^2(I)} = 1 + o(1)$ (thanks to similar computations as above), this proves the claimed lower bound.
\end{proof}

\subsection{Estimate for some pseudo-differential type operators on polynomials}\label{sec-pseudo-polyn-est}
In this section, we use the spectral analysis of the operator $P_\nu$ to deduce the operator estimate of \cref{th-pseudo-polyn-est}. In order to do that, we need some definitions and theorems about a general class of operators on polynomials. The following comes from~\cite[definition~9 \&\ theorem~18]{Koenig17}.

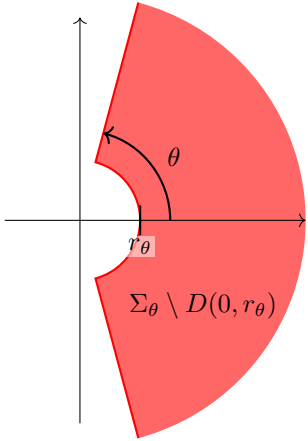
\begin{figure}
\begin{minipage}[c]{0.4\textwidth}
 \begin{center}
  \input{Delta_th.tex}
 \end{center}
 \end{minipage}\hfill%
 \begin{minipage}[c]{0.6\textwidth}
\caption{An example of a set $\sector[\th]\setminus D(0,r_\th) \subset \Omega$. The angle 
$\theta$ is allowed to be arbitrarily close to $\pi/2$, but then, the radius 
$r_\theta$ of the disk we avoid may blow up arbitrarily fast. For the $\Omega$ we will consider, the corresponding $r_\theta$ does blow up when $\theta \to \pi/2$ (at least, we cannot exclude that it blows up).}%TODO: ajouter Omega à la figure
\label{fig:U_r_theta}
\end{minipage}
\end{figure}

\begin{definition}\label{def-symbols}
Let $\Omega$ be an open subset of $\C$. Assume that there exists $(r_\th)_{0\leq\th<\pi/2}$ with $r_\th \geq 0$ such that $\bigcup_{0\leq\th<\pi/2}\sector[\th]\setminus D(0,r_\th) \subset \Omega$ (see \cref{fig:U_r_theta}).
  
We denote by $\Sc(\Omega)$ the set of functions $\gamma$ holomorphic on $\Omega$ that have sub-exponential growth on each $\sector[\th]\cap \Omega$, 
i.e., for each $\theta\in [0,\pi/2)$ and $\d>0$, we have 
\[
p_{\theta,\d}(\gamma) \coloneqq \sup_{\nu\in \sector[\th] \cap \Omega} |\gamma(z)\eu^{-\d|\nu|}| < +\infty.
\]
We endow $\Sc(\Omega)$ with the topology defined by the seminorms $p_{\theta,\d}$ 
for all $\theta\in[0,\pi/2)$ and $\d>0$.
\end{definition}

For the next theorem, if $U$ is an open subset of $\C$, we denote the set of bounded holomorphic functions on $U$ that have a zero of order $n_0$ at $0$ by $\Hol^\infty_{n_0}(U)$. We endow $\Hol^\infty_{n_0}(U)$ with the $L^\infty$-norm.
\begin{theorem}\label{th:est_default}
Let $\Omega\subset \C$ as in \cref{def-symbols} and set $n_0 = \min\{n\in\N\, : [n,+\infty) \subset \Omega\}$. Let $\gamma$ in $\Sc(\Omega)$ and $\Op{\gamma}$ be the operator on polynomials with a zero of order $n_0$ at $0$, defined by: 
\[\Op{\gamma}\bigg(\sum_{n\geq n_0} a_n z^n\bigg) = \sum_{n\geq n_0} \gamma(n) a_n z^n.\]
Let \(U\noic\) be a bounded open subset  of $\C$. Let \(V\noic\) be a neighborhood of \(\overline U\noic\) that is star shaped with respect to $0$. Then there exists $C>0$ such that for all polynomials $p$ with a zero of order $n_0$ at $0$,
\[\lVert\Op{\gamma}(p)\rVert_{L^\infty(U)}\leq C\|p\|_{L^\infty(V)}.\]

Moreover, the constant $C$ above can be chosen continuously in $\gamma\in \Sc(\Omega)$: the 
map $\gamma\in \Sc(\Omega) \mapsto \Op{\gamma}$ is continuous from $\Sc(\Omega)$ to $\Lc\left(\Hol^\infty_{n_0}(V),\, \Hol^\infty_{n_0}(U)\right)$.
\end{theorem}

We now have all the pieces needed to prove \cref{th-pseudo-polyn-est}.
\begin{proof}[Proof of \cref{th-pseudo-polyn-est}]
Let $\th_0 \in \big[0,\frac \pi 2 \big)$. According to \cref{prop-lambda-n} there exists $r_{\th_0}$ such that if $|\nu|>r_{\th_0}$ and $\lvert\arg(\nu)\rvert < \th_0$, then there exists a unique eigenvalue $\l_\nu$ of $P_\nu$ close to $q'(0)\nu$. Moreover, this eigenvalue is algebraically and geometrically simple.

Set $\Omega = \bigcup_{\th_0 \in [0,\pi/2)} \sector[\th_0] \setminus D(0,r_{\th_0})$. Notice that by definition, $\Omega$ satisfies the property of \cref{def-symbols}. For $0<t<T$, $x\in I$ and $\nu \in \Omega$, we define 
\[
\tilde\g_{t,x}(\nu) \coloneqq \eu^{-t(\l_\nu - q'(0)\nu)} \f_\nu(x) \eu^{\nu \dagmon(x)(1-\e)},
\]
and
\begin{equation}
    \label{eq-def-gamma-complex}
    \g_{t,x}(\nu) = \tilde \g_{t,x}(\nu+1)
\end{equation}
which is the natural extension of the definition of $\g_{t,x}(n)$ when $n\in\N$.
\step{The family $(\tilde \g_{t,x})_{0<t<T,x\in I}$ is a bounded family\footnote{Let us recall that if $E$ is locally convex vector space whose topology is generated by a family $(p_\iota)_\iota$ of seminorms, a subset $X$ of $E$ is bounded if and only if for every $\iota$, the set $\{p_\iota(x), x\in X\}$ is a bounded subset of $\R$.} of $\Sc(\Omega)$}
According to \cref{prop-lambda-n}, $\l_\nu$ is algebraically simple on $\Omega$. Thus, according to analytic perturbation theory (see, e.g., \cite[Chapter VII, \S1]{Kato95}), $\l_\nu$ and the associated spectral projection are holomorphic in $\nu \in \Omega$. Since $\f_\nu = \Pi_\nu \big(\nu^{1/4} \eu^{-\nu q'(0) x^2\!/2}\big)$, $\f_\nu$ is holomorphic in $\nu \in \Omega$. We conclude that $\tilde \g_{t,x}(\nu)$ is holomorphic in $\nu \in\Omega$.

We still have to prove that $(\tilde \g_{t,x})_{0<t<T,x\in I}$ is a bounded family of $\Sc(\Omega)$. Let us set 
\begin{gather*}
\eta_t(\nu) \coloneqq \eu^{-t(\lambda_\nu - q'(0)\nu)};\\
\zeta_x(\nu) \coloneqq \f_\nu(x) \eu^{\nu\dagmon(x)(1-\e)},
\end{gather*}
and prove that both of the families $(\eta_t)_{0<t<T}$ and $(\zeta_x)_{x\in I}$ are bounded in $\Sc(\Omega)$. It is easy to see that $(\g_1, \g_2)\in (\Sc(\Omega))^2 \mapsto \g_1\g_2 \in \Sc(\Omega)$ is bounded, so this will prove the claim.

Let $\th_0 \in \big[0,\frac \pi 2 \big)$ and $\d>0$. According to \cref{prop-lambda-n}, we have in the limit $|\nu|\to +\infty$, $\nu \in \sector[\th_0]$, $\l_\nu - \nu q'(0) = o(\nu)$. Thus, for $\nu$ large enough in $\sector$,
\[
 |\eta_t(\nu)\eu^{-\d|\nu|}| = |\eu^{ t o(|\nu|) - \d |\nu|}| < \eu^{TC_{\th_0,\d}}.
\]
Thus, $(\eta_t)_{0<t<T}$ is a bounded family of $\Sc(\Omega)$.

Similarly, according to \cref{prop-agmon-linfty}, we have for any $x\in I$, and $\nu$ large enough in $\sector$,
\[
|\zeta_x(\nu)\eu^{-\d|\nu|}| \leq C.
\]
This prove that $(\zeta_x)_{x\in I}$ is a bounded family of $\Sc(\Omega)$.

\step{The family $(\g_{t,x})_{0<t<T,x\in I}$ is a bounded family of $\Sc(\Omega)$}
According to the definition of $\Omega$ as a union of domains that look like the one of \cref{fig:U_r_theta}, $\Omega$ is stable by $\nu\mapsto \nu+1$. Then, the map $\g \in \Sc \mapsto \g(\cdot+1) \in \Sc$ is well-defined and continuous. Thus, according to the first step and the definition of $\g$ (\cref{eq-def-gamma-complex}), the family $(\g_{t,x})_{0<t<T,x\in I}$ is indeed a bounded family of $\Sc(\Omega)$.

\step{Conclusion}
Let $n_0 = \min\{n\in\N\, : [n,+\infty) \subset \Omega\}$.\footnote{In fact, we can be more precise in the the construction of $\Omega$ and ensure that $\R_+\subset \Omega$, in which case $n_0 = 0$. Indeed, the spectral theory of compact operators proves that for every $\nu\in \C$, the spectrum of $P_\nu$ is a discrete sequence of eigenvalues. The uniqueness of the solution of Cauchy problems for ODEs proves that when $\nu>0$, these eigenvalues are actually geometrically and algebraically simple. Finally, perturbation theory proves that the first eigenvalue is holomorphic on the neighborhood of $\R_+$. We do not need this, so we do not detail this.} 
Let $U$ be a bounded open neighborhood of $X$ such that $\overline U\subset V$. Then, the sets $U$ and $V$ satisfy the hypotheses of \cref{th:est_default}. Hence, according to \cref{th:est_default}, there exists $C>0$ such that for every polynomials $p$ with a zero of order $n_0$ at $0$, and for every $x\in I$ and $0<t<T$,
\[
\lVert\Op{\g_{t,x}}(p)\rVert_{L^\infty(U)}\leq C\|p\|_{L^\infty(V)}. \qedhere
\]
\end{proof}

\section{Critical time of null-controllability for some domains}\label{sec-crit-time}

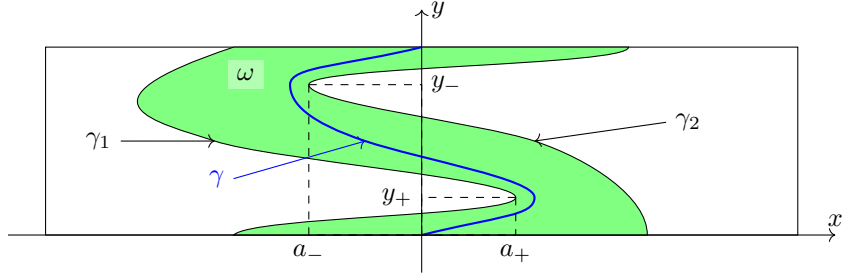
\begin{figure}
\centering
\input{omega_corollary}
\caption{In green, the domain $\omega$. At $y=y_- $, the function $\gamma_2^-$ takes its maximum $a_-$. Then, the interval $(a_-,L_+)\times \{y_-\}$ is disjoint from $\overline{\omega}$. So, the Grushin equation is not null-controllable in time $T<\dagmon(a_-)/q'(0)$. Similarly, the interval $(-L_-, a_+)\times \{y_+\}$ is disjoint from $\overline{\omega}$. So, the Grushin equation is not null-controllable in time $T<\dagmon(a_+)/q'(0)$. Also, if we take a path $\gamma$ (here in blue) that  is close to the boundary of $\omega$ around $y = y_-$ and $y = y_+$, then, we can apply \cref{th-main-pos-intro}, and the Grushin equation is null-controllable in time $T> \max(\dagmon(a_-),\dagmon(a_+))/q'(0)$.}
\label{fig:omega_corollary}
\end{figure}

In this section, we prove \cref{th-main}. 

\begin{proof}[Proof of \cref{th-main}]
Set $a_- = -\max(\g_2^-)$ and $a_+ = \max(\g_1^+)$. Denote by $y_-\in \T$ and $y_+\in \T$ points where these maxima are reached.
\step{Lower bound of the minimal time} For this step, we only have to treat the case $T_*>0$. In this case, either $a_-< 0$ or $a_+ >0$. If $a_-<0$, for any $a_- < a < 0$, the segment $[a,L_+]\times \{y_-\}$ stays at
positive distance of $\omega$, and thanks to \cref{th-main-neg}, the generalized Baouendi-Grushin equation~\eqref{eq-grushin-gen} is not null-controllable on $\omega$ in time $T<\dagmon(a)/q'(0)$. 
Similarily, if $a_+ >0$, for any $0<a<a_+$, the segment $[-L_-,a]\times\{y_+\}$ stays at positive distance from $\o$, and the generalized Baouendi-Grushin equation is not null-controllable in time $T<\dagmon(a)/q'(0)$.

This holds for any $a_-<a<0$ and $0<a<a_+$, thus the generalized Baouendi-Grushin equation~\eqref{eq-grushin-gen} is not null-controllable in time $T<T_*$.

\step{Upper bound of the minimal time} Let $\varepsilon>0$ small enough so that $\gamma_2-\gamma_1 > \varepsilon$. Let $\tilde\gamma_1 =\max(\gamma_1, a_--\varepsilon)$ and $\tilde\gamma_2 =\min(\gamma_2, a_++\varepsilon)$. By using the information $\gamma_2-\gamma_1>\varepsilon$, $a_- \leq \gamma_2$, $\gamma_1\leq a_+$ and by looking at the different cases, we readily get $\tilde \gamma_2 - \tilde \gamma_1 \geq \varepsilon$. Then, we define the path 
\[
\gamma = (\g_x,\g_y)\colon s\in \T \mapsto \left(\frac{\tilde\gamma_1(s) + \tilde\gamma_2(s)}2, s\right).
\]
This path satisfies $\gamma_1 +\varepsilon/2 \leq \gamma_x \leq \gamma_2-\varepsilon/2$, hence $\gamma(\T) \subset  \o$. Moreover, we see that it satisfies the hypotheses of \cref{th-main-pos-intro}, because the connected components of $(I\times \T) \setminus \g(\T)$ are $\{(x,y)\colon x< \g_x(y)\}$ and $\{(x,y) \colon x>\g_x(y)\}$. Moreover, this path satisfies 
\[
a_- -\frac\e2 \leq \gamma_x \leq a_+ +\frac\e2. 
\]
Thus,  \cref{th-main-pos-intro} implies that the generalized Baouendi-Grushin equation~\eqref{eq-grushin-gen} is null-controllable in time $T>\max(\dagmon(a_++\e/2),\dagmon(a_--\e/2))/q'(0)$. As this holds for every $\e>0$ small enough, the result follows.
\end{proof}

\appendix

\section{Control of the Baouendi-Grushin equation on $I\times(0,\pi)$}\label{app-dirichlet}
In the article, we stated and proved results on the Baouendi-Grushin posed on $I\times \T$. These results have a version for the Baouendi-Grushin equation posed on $I\times (0,\pi)$:
\begin{equation}
 \label{eq-grushin-gen-dirichlet}
 \left\{\begin{array}{ll}
         (\partial_t-\partial_x^2 - q(x)^2 \partial_y^2)f(t,x,y) = \1_\o u(t,x,y),&\quad t\in (0,T), x\in I, y\in (0,\pi)\\
         f(t,x,y) = 0,&\quad t\in (0,T), (x,y)  \in \partial (I\times (0,\pi)), \\
         f(0,x,y) = f_0, & \quad x\in I, y\in (0,\pi).
        \end{array}\right.
\end{equation}
Here, we precisely state them and explain what are the differences, if any, in their proofs. The precise definition of the operator, especially its domain, is again the Friedrichs' extension. That it generates an analytic semigroup is again proved with Hille-Yosida's theorem. We again refer to~\cite{Helffer13,Brezis11,BCG14} for the details.

The adaptation of the positive controllability result is:
\begin{theorem} \label{th-main-pos-dirichlet}
	Assume that $q\in C^3(\overline I)$ is such that $q(0) = 0$ and $\min_{\overline I} q' >0$. Let $\o$ be an open subset of $I\times [0,\pi]$. Assume that there exists $\gamma = (\g_x, \g_y) \in C^0([0,1], I\times [0,\pi])$ such that $\g((0,1)) \subset \o$, $\g_y(0) = 0$ and $\g_y(1) = \pi$.
	
	The generalized Baouendi-Grushin equation~\eqref{eq-grushin-gen-dirichlet} is null-controllable on $\o$ in time $T$ such that 
	\[
	T> \frac{1}{q'(0)}\max\left(\dagmon\left(\min_\T(\g_x)\right), \dagmon\left( \max_\T (\g_x)\right)\right).
	\]
\end{theorem}
The proof is mostly the same, the only small difference being the construction of the cutoff function, which is done thanks to \cref{th-component-segment} and the natural adaptation of \cref{th-cutoff}.

The adaptation of the negative result~\cref{th-main-neg-intro} is straightforward. We use the same notation $\d$ as in \cref{th-main-neg-intro}:
\begin{theorem} \label{th-main-neg-dirichlet}
	Assume that $q\in C^2(\overline I)$ is such that $q(0) = 0$, $q'(0) > 0$ and $q(x) \neq 0$ whenever $x\neq 0$.
	Let $\o$ be an open subset of $I\times (0,\pi)$. Assume that there exist $a \in [-L_-,0)$, $b \in (0,L_+]$ and $y_0 \in (0,\pi)$ such that 
	\[
	\distance\big((a,b)\times \{y_0\}, \o\big) > 0.
	\]
	Then, the generalized Baouendi-Grushin equation~\eqref{eq-grushin-gen-dirichlet} is not null-controllable on $\o$ in time $T$ such that 
	\[
	T< \frac{1}{q'(0)}\min \left(\dagmontilde(a),\dagmontilde(b) \right).
	\]
\end{theorem}
We prove this theorem with the following lemma:
\begin{lemma}
Let $T>0$ and let $\omega$ be an open subset of $I\times (0,\pi)$. Denote by $S(\omega)$ the symmetric of $\omega$ with respect to $\{y = 0\}$. Let $\l_n$ and $\f_n$ as in \cref{eq-def-vn}.

Assume that the generalized Baouendi-Grushin equation with Dirichlet boundary conditions~\eqref{eq-grushin-gen-dirichlet} is null controllable on $\omega$ in time $T$, then for every complex sequence $(a_n)$ with a finite number of nonzero terms,
\[
\sum_{n\geq 0} |a_n|^2 \|\f_n\|^2_{L^2(I)}\eu^{-2\l_n T} \leq C \int_{[0,T]\times (\o \cup S(\o))} \bigg|\sum_{n\geq 0} a_n \f_n(x) \eu^{\iu n y - \l_n t}\bigg|^2 \diff t \diff x \diff y.
\]
\end{lemma}
\begin{proof}[Sketch of the proof]
This lemma is proved by testing the associated observability inequality on the function $g(t,x,y) = \sum_{n\geq 0} a_n \f_n(x)\sin(ny) \eu^{-\l_n t}$, and writing $\sin(ny) = (\eu^{\iu ny} - \eu^{-\iu ny})/(2\iu)$. Thus, with $\tilde g(t,x,y) = \sum_{n\geq 0} a_n \f_n(x)\eu^{\iu n y - \l_n t}$, $g(t,x,y) = (\tilde g(t,x,y) - \tilde g(t,x,-y))/(2\iu)$, the right-hand side of the observability inequality satisfies
\[
\|g\|_{L^2([0,T]\times \o)}^2 \leq \dfrac12(\|\tilde g\|_{L^2([0,T]\times \o)}^2 + \|\tilde g\|^2_{L^2([0,T]\times S(\o))}).
\]
The right-hand side of this inequality is the right-hand side of the claimed estimate.
\end{proof}
\Cref{th-main-neg-dirichlet} is then proved by remarking that we already disproved such an inequality in \cref{sec-sketch-gen}.

\section{Existence of the cutoff function and homotopy}\label{app-cutoff}
We begin with the construction of the cutoff function used in the proof of \cref{th-main-pos-intro}.
\begin{proposition}\label{th-cutoff}
Let $a<b$ and let $\o$ be an open subset of $(a,b)\times \T$. Assume that there exists a closed path $\gamma = (\g_x, \g_y) \in C^0(\T;\,\o)$ such that $\{a\}\times \T$ and $\{b\}\times \T$ are included in different connected components of $([a,b]\times\T) \setminus \g(\T)$. Let $\o_- = [a,\min \g_x]\times \T$ and $\o_+ = [\max \g_x,b]\times \T$. There exists a function $\chi\in C^\infty([a,b]\times \T)$ such that:
\begin{itemize}
    \item $\chi = 0$ on $\o_+ \setminus \o$;
    \item $\chi = 1$ on $\o_- \setminus \o$;
    \item $\supp(\nabla \chi)\subset \o$.
\end{itemize}
\end{proposition}

\begin{proof}[Proof of \cref{th-cutoff}]
\step{Defining $\chi$} Let $\epsilon>0$ small enough so that $\distance(\g(\T),\{a,b\}\times \T)>\epsilon$ and such that for any $z\in \g(\T)$, $B(z,\epsilon) \subset \o$. We set
\[
\o_0 \coloneqq \{z\in[a,b] \times \T\colon \distance(z,\g(\T)) < \epsilon\}.
\]
Let $\rho \in C_c^\infty(B(0,\epsilon/2))$ with $\int_{B(0,\epsilon/2)} \rho(z)\diff z = 1$. Consider $\O$ the connected component of $\{a\}\times\T$ in $(\R\times \T)\setminus \g(\T)$. Set $\chi = \rho \ast \1_{\O}$ (initially defined on $\R\times \T$ and then restricted on $[a,b]\times \T$).

\step{$\supp(\nabla \chi) \subset \o_0$} According to the definition of $\O$, $\1_\O$ is locally constant outside of $\g(\T)$. This implies that $\chi$ is locally constant around each $z$ such that $\distance(z,\g(\T))>\epsilon/2$. Hence $\supp(\nabla \chi)\subset \{z\colon \distance(z,\g(\T))\leq \epsilon/2\}$. According to our choice of $\epsilon$, this proves the claim that $\supp(\nabla \chi)\subset \o_0 \subset \o$.

\step{Value of $\chi$ on $\o_-\setminus \o$} According to the definition of $\chi$ and the fact that $\distance(\{a\}\times\T,\g(\T))> \epsilon$, for any $y_0\in\T$, $\chi(a,y_0) = 1$. Moreover, $\o_-\setminus \o_0$ is connected (according to the definition of $\o_-$, we can connect every $(x,y) \in \o_-\setminus \o_0$ to the left boundary $\{a\}\times \T$ with the horizontal segment $[a,x]\times\{y\}$). According to the previous step, $\chi$ is locally constant outside $\o_0$. Hence, $\chi$ is constant in $\o_-\setminus \o_0$.

\step{Value of $\chi$ on $\o_+\setminus \o$} According to the definition of $\chi$ and \cref{th-component}, $\chi = 0$ on $\{b\}\times \T$. The rest of this step is a copy-paste of the previous step.
\end{proof}

Now, we justify \cref{rk-homotopy}, with the following two propositions:
\begin{proposition}\label{th-component}
Let $a<b$ and let $\g \in C^0(\T, (a,b)\times \T)$ be a closed path that is not homotopic to a constant path. Then $\{a\}\times \T$ and $\{b\}\times \T$ are included in different connected components of $([a,b]\times\T) \setminus \g(\T)$.
\end{proposition}

\begin{proof}
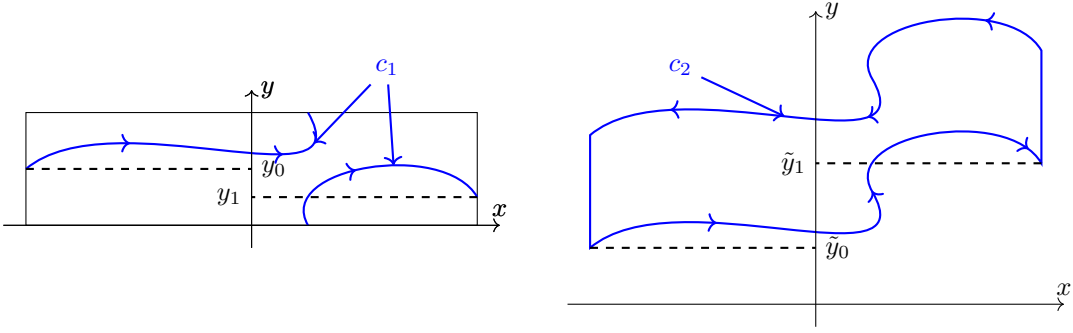
\begin{figure}
\centering
\input{covering}
\caption{Illustration of the path $c_2$ defined in the proof of \cref{th-component}.}
\label{fig:covering}
\end{figure}
Assume that for some $y_0, y_1\in\T$, there exists a continuous path $c_1$ in $(\overline{I} \times \T) \setminus \g(\T)$ from $(a,y_0)$ to $(b,y_1)$. Since $\overline I\times \T$ is Hausdorff, we may assume that $c_1$ is simple. We can also assume that $c_1$ touches the boundaries $\{a,b\}\times\T$ only at the start and end.

Now, consider the universal cover $[a,b]\times \R$ of $[a,b]\times \T$. Consider $\tilde c_1$ a lift of $c_1$ to $[a,b]\times \T$, that starts at $(a,\tilde y_0)$ and ends at $(b,\tilde y_1)$. Let $c_2$ the simple closed loop formed by concatenating $c_1$, the vertical segment $\{b\}\times [\tilde y_1,\tilde y_1+2\pi]$, the reverse of the path $\tilde c_1 +(0,2\pi)$, and finally the vertical segment $\{a\}\times [\tilde y_0,\tilde y_0+2\pi]$ from top to bottom (see \cref{fig:covering}).

If we see this path $c_2$ as a path on $\R^2$, according to Jordan's theorem, $\R^2\setminus c_2$ has two path-connected components, one of them bounded. Let us denote by $\O_1$ this bounded component, which, according to Jordan-Schoenflies' theorem, is simply connected. One of the lift of $\g$ lies in $\O_1$, let us call it $\tilde \g$. But $\g$ is not homotopic to a constant path, which contradicts the simple connectedness of $\O_1$.
\end{proof}

\begin{proposition}\label{th-homotopy-2}
Let $\o$ be a connected open subset of $[a,b]\times \T$ such that $\{a\}\times \T$ and $\{b\}\times \T$ are included in different connected components of $([a,b]\times\T) \setminus \o$. Let $\widetilde \o$ a connected open subset of $[a,b]\times \T$ such that $\overline{\o}\subset \widetilde \o$.

There exists a closed path $\g \in C^0(\T,\widetilde \o)$ that is not homotopic in $[a,b]\times \T$ to a constant path.
\end{proposition}
\begin{remark}\label{rk-homotopy-2}
Let $\g$ be a closed path in $(a,b) \times \T$ such that $\{a\}\times \T$ and $\{b\}\times \T$ are included in different connected components of $([a,b]\times\T) \setminus \g(\T)$. It is possible this path $\g$ is homotopic to a constant path, but \cref{th-homotopy-2} applied with $\o \coloneqq \{z\colon \distance(z,\g(\T)) < \e\}$ and $\tilde \o \coloneqq \{z\colon \distance(z,\g(\T)) < 2\e\}$ tells us that for any $\e > 0$ there exists a path $\tilde \g$ that stays at distance at most $2\e$ from $\g(\T)$ and that is not homotopic to a constant path.
\end{remark}
\begin{proof}
The proof uses some basic tools of algebraic topology, in particular van Kampen's theorem (see for instance Hatcher's ``Algebraic Topology''~\cite[\S 0.1, \S1.1--1.2]{Hatcher02}).

Let $C_-$ be the connected component of $\{a\}\times \T$ in $([a,b]\times \T) \setminus \o$. Set 
\begin{equation}\label{eq-defA+A-}
    A_- \coloneqq C_- \cup \widetilde \o, \qquad A_+ \coloneqq \big(([a,b]\times \T)\setminus C_-\big)\cup \widetilde \omega.
\end{equation}

\step{Every connected component $C$ of $([a,b]\times \T) \setminus \o$ is closed and satisfies $\partial C \subset \tilde \o$}\label{step-closed} Here, $\partial C$ is the boundary of $C$ as a subset of $[a,b]\times \T$.

Indeed, connected components of a topological space $X$ are closed in $X$, hence $C$ is closed in $([a,b]\times \T) \setminus \o$, which is itself closed in $[a,b]\times \T$.

If $x\in \partial C$ is such that $x\notin \partial \o$, then there exists $\e>0$ such that $B(x,\e)\subset ([a,b]\times \T) \setminus \o$. But then, $B(x,\e)$ is included in the connected component of $x$, i.e., $C$. By contradiction, we see that every $x\in \partial X$ is in $\partial \o \subset \tilde \o$.

\step{$A_+$ and $A_-$ are open and connected}
We begin with the openness of $A_-$. If $x \in A_-$, there are three cases:
\begin{itemize}
    \item If $x$ is in the interior of $C_-$, it is in the interior of $A_-$ by definition.
    \item If $x$ is in $\widetilde \o$, since $\widetilde \o$ is open, $x$ is also in the interior of $A_-$.
    \item If $x\in \partial C_-$, according to step \ref{step-closed}, $x\in \tilde \o$ which implies that $x$ is in the interior of $A_-$. 
\end{itemize}
The subset $A_+$ is open because it is the union of the open subsets $([a,b]\times \T)\setminus C_-$ and $\tilde \o$.

Since $A_-$ is the union of two connected subset that have a non-empty intersection (we saw in step \ref{step-closed} that $\partial C_- \subset C_-\cap \widetilde\o$), $A_-$ is connected. Finally, $A_+$ is connected because it is the union of the connected subset $\widetilde \o$ and of the connected components of $[a,b]\times \T \setminus \o$ other than $C_-$, which all have a non-empty intersection with $\widetilde \o$.

\step{Conclusion using van Kampen's theorem}
If $\alpha$ is a closed path in a topological space $X$, we will denote its homotopy class by $[\a]_X$. We will denote the fundamental group of $X$ by $\pi_1(X)$. We will denote by $p_-$ (respectively $p_+$) the canonical injection of $\pi_1(A_-)$ (respectively $\pi_1(A_+)$) into the free product $\pi_1(A_-) \ast \pi_1(A_+)$.

According to van Kampen's theorem~\cite[Theorem~1.20]{Hatcher02}, the map
$k\colon\pi_1(A_-)\ast \pi_1(A_+) \to \pi_1([a,b]\times \T)$ defined by 
\[
[\a]_{A_\pm} \in \pi_1(A_\pm) \mapsto [\a]_{[a,b]\times\T} \in \pi_1([a,b]\times \T)
\]
is surjective. Moreover, its kernel is generated by $p_+([\a]_{A_+}) p_-([\a]_{A_-})^{-1}$ for all closed paths $\a$ in $A_-\cap A_+ = \widetilde \o$.

If we denote the closed path $s\in\T \mapsto (a,s)\in [a,b]\times\T$ (respectively $s\mapsto (b,s)$) by $\b_-$ (respectively $\b_+$), the definition of $k$ implies that
\[
k(p_-([\b_-]_{A_-})) = [\b_-]_{[a,b]\times \T} = [\b_+]_{[a,b]\times \T} = k(p_+([\b_+]_{A_+})).
\]
Thus, $\xi \coloneqq p_-([\b_-]_{A_-}) p_+([\b_+]_{A_+})^{-1} \in \ker(k)$. According to the previous discussion, $\xi$ is a product of terms of the form $p_\pm([\a_k]_{A_\pm}) p_\mp([\a_k]_{A_\mp})^{-1}$ for a finite number of paths $\a_k$ in $\widetilde \o$. Reducing the word $\xi$ in the free product $\pi_1(A_-)\ast \pi_1(A_+)$, $\xi$ can be written in the form
\begin{equation} \label{eq:xi}
p_-([\b_-]_{A_-}) p_+([\b_+]_{A_+})^{-1} = \xi = p_\pm([\tilde \a_1]_{A_\pm}) p_\mp([\tilde \a_2]_{A_\mp})p_\pm([\tilde \a_3]_{A_\pm})\cdots
\end{equation}
where none of the terms in the right-hand side are the neutral element of $\pi_1(A_\pm)$. Since the left-hand side is already a reduced word, by definition of the free product of groups, the two words on the left and right-hand side of this equality are the same. Thus there are exactly two factors in the right-hand side of \cref{eq:xi} and
\begin{gather*}
    [\b_-]_{A_-} =[\tilde \a_1]_{A_-},\quad [\b_+]_{A_+}^{-1}= [\tilde \a_2]_{A_+}.
\end{gather*}
The first equality tells us that $\tilde \a_1$ is homotopic in $A_-$ to $\b_-$. Since $\tilde \a_1$ is in $\widetilde \o$ and since $\b_-$ is not homotopic to a constant path in $[a,b]\times \T$, this proves the proposition.
\end{proof}

To end this appendix, we mention that \cref{th-component} has a variant when the domain is $(a,b)\times [0,\pi]$ instead of $(a,b)\times \T$:
\begin{proposition}\label{th-component-segment}
Let $a<b$ and let $\g = (\g_x,\g_y) \in C^0([0,1], (a,b)\times [0,\pi])$ be a closed path such that $\g_y(0) = 0$ and $\g_y(1) = \pi$. Then $\{a\}\times [0,\pi]$ and $\{b\}\times [0,\pi]$ are included in different connected components of $[a,b]\times[0,\pi] \setminus (\g([0,1]))$.
\end{proposition}
The proof also uses Jordan's theorem, but in a simpler way than \cref{th-component}, and is left to the reader.

\section{Non-selfadjoint harmonic oscillators} \label{sec:Davies-operator}

Let $\b \in \C$ with $\Re(\b) > 0$. We discuss in this appendix the basic properties of the non-selfadjoint harmonic oscillator (or Davies operator) defined on $L^2(\R)$ by
\begin{equation} \label{def:Hbeta}
H_\beta = - \partial_x^2 + \beta^2 x^2.
\end{equation}
More precisely, we set 
\[
\Dom(H_\beta) = \set{u \in L^2(\R) \colon  {-}\partial_x^2 u + \beta^2 x^2 u \in L^2(\R)},
\]
where $-\partial_x^2 u + \beta^2 x^2 u$ is understood in the sense of distributions, and we define $H_\b$ by \cref{def:Hbeta} on $\Dom(H_\beta)$. This defines an unbounded operator on $L^2(\R)$. When $\b = 1$ we recover the usual harmonic oscillator. 

The spectral properties of the Davies operator has been studied (see, among others, \cite[\S 14.5]{Davies07}, \cite[\S 14.4]{Helffer13}, \cite{KSTV15}, and the references therein), and the properties stated in this appendix are standard, at least in spirit. Nevertheless, for the reader convenience, we collect and prove the properties needed in our study.

\begin{proposition} \label{prop:Davies}
Let $\b \in \C$ with $\Re(\b) > 0$.
\begin{enumerate}[\rm (i)]
\item The operator $H_\b$ is closed and has compact resolvent.

\item The adjoint of $H_\b$ is $H_\b^* = H_{\bar \b}$.

\item The spectrum of $H_\beta$ is given by the sequence of (geometrically and algebraically) simple eigenvalues $(\l_{\b,k}^\h)_{k \in \N^*}$ where 
\[
\l_{\b,k}^\h = (2k-1) \beta.
\]
An eigenfunction associated to $\l_{\b,k}^\h$ is given by 
\[
\f_{\b,k}^\h(x) = \He_{k-1}(\sqrt \b x)\eu^{-\b x^2\!/2}, \quad \text{where} \quad \He_{k-1}(x) = (-1)^{k-1}\eu^{x^2}\partial_x^{k-1}(\eu^{-x^2})
\]
is the $(k-1)$-th Hermite polynomial. 

\item 
For $\g > 0$ and $\e > 0$ we set 
\begin{equation} \label{def:Zc-beta}
\Zc_{\b,\e,\g} \coloneqq \{ \nu \in \C \colon \abs z \leq \g \abs \beta, \, \distance (z ,\sigma(H_\beta)) \geq \e \abs \nu\}.
\end{equation}
Let $\th_0 \in \big[0,\frac \pi 2 \big)$. There exists $C > 0$ such that if $\beta \in \sector$ (see \cref{def:sector}) then for $z \in \Zc_{\b,\e,\g}$ we have 
\begin{equation} \label{eq:res-Hb}
\big\| (H_\beta - z)\inv \big\|_{\Lc(L^2(\R))} \leq \frac C {\abs \beta}.
\end{equation}

\item We denote by $\Pi_{\b}^\h$ the spectral projection of $H_\b$ associated to the eigenvalue $\l_{\b,1}^\h = \b$. Then we have
\begin{equation}\label{eq:eigenproj-Hb-formula}
    \Pi_\b u = \frac {\big\langle \f_{\bar \b,1}^\h,u \big\rangle} {\big\langle \f_{\bar \b,1}^\h , \f_{\b,1}^\h \big\rangle} \f_{\b,1}^\h.
\end{equation}
and
\begin{equation} \label{eq:eigenproj-Hb}
\nr{\Pi_\b^\h}_{\Lc(L^2(\R))} = \sqrt{\frac {\abs \b}{\Re(\b)}} = \frac 1 {\sqrt{\cos(\arg(\b))}}.
\end{equation}
\end{enumerate}
\end{proposition}

\begin{proof}
\step{$H_\b$ has compact resolvent} The closedness of $H_\b$ is clear. Let $\th = \arg(\beta) \in \big( - \frac \pi 2,\frac \pi 2 \big)$. For $u \in \Dom(H_\b)$ we have 
\begin{equation} \label{eq:innp-Hb}
\innp{\eu^{-\iu\th} H_\b u}{u}_{L^2(\R)} = \eu^{-\iu\th} \nr{u'}_{L^2(\R)}^2 + \abs \beta^2 \eu^{\iu\th} \nr{xu}_{L^2(\R)}^2,
\end{equation}
so $\eu^{-\iu\th} H_\b$ is sectorial with angle $\th$. In particular, $(\eu^{-\iu\th} H_\b + 1)$ is injective. Now let
\[
B^1(\R) = \set{u \in H^1(\R) \colon  xu  \in L^2(\R)}.
\]
This is a Hilbert space for the natural norm
\[
\nr{u}_{B^1(\R)}^2 = \nr{u}_{H^1(\R)}^2 + \nr{xu}_{L^2(\R)}^2.
\]
For $u,v \in B^1(\R)$ we set 
\[
Q_\beta(u,v) = \eu^{-\iu\th} \innp{u'}{v'}_{L^2(\R)} + \abs \beta^2 \eu^{\iu\th} \innp{xu}{xv}_{L^2(\R)} + \innp{u}{v}_{L^2(\R)}.
\]
Let $f \in L^2(\R)$. By the Lax-Milgram Theorem, there exists a unique $u \in B^1(\R)$ such that $Q_\b(u,v) = \innp{f}{v}_{L^2(\R)}$ for all $v \in B^1(\R)$. In the sense of distributions we have $\eu^{-\iu\th} (-u'' + \b^2 x^2 u) + u = f \in L^2(\R)$, so $u \in \Dom(H_\b)$ and $(\eu^{-\iu\th} H_\b + 1)u = f$. This proves that $\eu^{\iu\th}$ belongs to the resolvent set of $H_\b$.

Finally, taking the real part of \cref{eq:innp-Hb} gives for $u \in \Dom(H_\b)$
\[
\nr{u'}_{L^2(\R)}^2 + \abs {\b}^2 \nr{xu}_{L^2(\R)}^2 \leq \frac {\Re \innp{\eu^{-\iu\th} H_\b u}{u}}{\cos(\th)} \leq \frac 1 {2 \cos(\th)} \left( \nr{H_\b u}_{L^2(\R)}^2 + \nr{u}_{L^2(\R)}^2  \right) .
\]
We deduce that $\Dom(H_\b)$ is compactly embedded in $L^2(\R)$. Since $H_\b$ has nonempty resolvent set, it has compact resolvent. In particular, its spectrum consists of isolated eigenvalues of finite multiplicities.

\step{Computation of $(H_\b)^*$} Let $v \in \Dom(H_{\bar \b})$. For $u \in \Dom(H_\b)$ we have $\innp{H_\b u}{v} = \innp{u}{H_{\bar \b}v}$, so $\Dom(H_{\bar \b}) \subset \Dom(H_\b^*)$ and $H_\b^* = H_{\bar \b}$ on $\Dom(H_{\bar \b})$. Now let $v \in \Dom(H_\b^*)$ and $f = H_\b^* v$. Then $v \in L^2(\R)$ and in the sense of distributions we have $-v'' + \bar \b^2 x^2 v = f \in L^2(\R)$, so $v \in \Dom(H_{\bar \b})$. This proves that $\Dom(H_\b^*) \subset \Dom(H_{\bar \b})$, and hence $H_\b^* = H_{\bar \b}$.

\step{Eigenvalues and eigenfunctions of $H_\b$, completeness of the eigenfunctions}\label{step-eigen-harmonic} For $k \in \N^*$ we have $\f_{\b,k}^\h \in \Dom(H_\b)$ and it is classical computation that for $\beta = 1$, (see, e.g., \cite[\S1.3]{Helffer13})
\[
H_1 \f_{1,k}^\h = \l_{1,k}^\h \f_{1,k}^\h
\]
Noticing that $\f_{\b,k}^\h(x) = \f_{1,k}(\sqrt\b x)$, routine computations using the scaling $x' = \sqrt \b x$ show that
\[
H_\b \f_{\b,k}^\h = \l_{\b,k}^\h \f_{\b,k}^\h.
\]
Then $\l_{\b,k}^\h$ is an eigenvalue of $H_\b$ and $\f_{\b,k}^\h$ is a corresponding eigenfunction. 

Let $u \in \Span(\f_{\b,k}^\h)_{k \in \N^*}^\bot$. Then for all polynomial $p$ we have $\int_\R u(x) p(\sqrt \b x) \eu^{-\b x^2\!/2} \diff x = 0$. For $\xi \in \R$ we set $F(\xi) = \int_\R \eu^{\iu x\xi} u(x) \eu^{-\b x^2\!/2} \diff x$. Then $F$ is analytic and $F^{(n)}(0) = 0$ for all $n \in \N$. This implies that $u(x) = 0$ for almost all $x \in \R$, so the family $(\f_{\b,k}^\h)_{k \in \N^*}$ is complete.

\step{Resolvent estimate} The map 
\[
(\th, \z) \mapsto \big(H_{\eu^{\iu\th}} - \eu^{\iu\th} \z \big)\inv
\]
is continuous and hence bounded on the compact $[-\th_0,\th_0] \times \Zc_{1,\e,\g}$, so \cref{eq:res-Hb} holds if $\abs \b = 1$.

For $\rho > 0$ we consider on $L^2(\R)$ the unitary operator $\Th_\rho$ such that for $u \in L^2(\R)$ and $x \in \R$ we have 
\begin{equation*} 
(\Th_\rho u)(x) = \rho^{\frac 12} u (\rho x).
\end{equation*}
We observe that  
\[
\Th_{\abs \beta^{\frac 12}}\inv H_\beta \Th_{\abs \beta^{\frac 12}} = \abs \beta H_{\frac \b {\abs \b}},
\]
so for $z \in \Zc_{\b,\e,\g}$ we have  
\[
\nr{(H_\b-z)^{-1}}_{\Lc(L^2(\R))} = \frac 1 {\abs \beta} \nr{\left( H_{\frac \b {\abs \b}}  - \frac z {\abs \b} \right)^{-1}}_{\Lc(L^2(\R))}.
\]
Since $\abs \b^{-1} \Zc_{\b,\e,\g}= \Zc_{\b/\abs \b,\e,\g}$, we deduce \cref{eq:res-Hb} in the general case.

\step{Spectral projection}
In the integral $\big\langle\f^\h_{\overline{\b},k},\f^\h_{\b,1}\big\rangle = \int_\R \He_{k-1}(\sqrt \b x)\eu^{-\b x^2}\diff x$, we make the change of variables and integration path $x' = \sqrt{\b} x$, which can be justified thanks to the gaussian decay of the integrand, and we find
\[
 \big\langle\f^\h_{\overline{\b},k},\f^\h_{\b,1}\big\rangle
 = \frac1{\sqrt \b} \int_\R \He_{k-1}(x) \eu^{-x^2}\diff x
 = \frac1{\sqrt \b} \innp{\f^\h_{1,k}}{\f^\h_{1,1}}.
\]
Since the functions $(\f^\h_{1,k})$ are eigenfunctions of the self-adjoint operator $H_1$ associated to different eigenvalues, we have $\langle \f^\h_{1,k},\f^\h_{1,1} \rangle = 0$ for $k\neq 1$. Since $\langle \f^\h_{1,1},\f^\h_{1,1} \rangle = \sqrt{\pi}$, we finally have
\begin{equation} \label{eq:phi-phi-bar}
 \big\langle\f^\h_{\overline{\b},k},\f^\h_{\b,1}\big\rangle
 =\left\{\begin{aligned}
           \sqrt{\frac\pi\b}&\ \text{ if } k = 1;\\
           0 &\ \text{ if } k \neq 1.
          \end{aligned}\right.
\end{equation}
Thus, we can define
\[
\widetilde \Pi_\b u = \frac {\big\langle \f_{\bar \b,1}^\h,u \big\rangle} {\big\langle \f_{\bar \b,1}^\h , \f_{\b,1}^\h \big\rangle} \f_{\b,1}^\h.
\]
Then, for $k\in\N^*$
\[
\widetilde \Pi_{\b} \f_{\b,k}^\h = \left\{\begin{array}{l}
\f_{\b,1}^\h \text{ if } k = 1;\\ 0 \text{ if } k \neq 1.
\end{array}\right.
\]
According to step \ref{step-eigen-harmonic}, the family $(\f_{\b,k}^\h)_k$ is complete, hence, by density, $\widetilde \Pi_\b$ is indeed the spectral projection $\Pi^\h_\b$.

Since $\f_{\b,1}^\h(x) = \eu^{-\b x^2\!/2}$, we can compute 
\[
\|\f_{\b,1}^\h\|_{L^2(\R)}^2 = \|\f_{\bar \b,1}^\h\|_{L^2(\R)}^2 = \sqrt{\frac {\pi}{\Re(\b)}}
\]
and \cref{eq:eigenproj-Hb} follows with \cref{eq:phi-phi-bar}.
\end{proof}

\begin{corollary}\label{th-trace-eigenproj}
Let $\Pi_\b^\h$ as in \cref{prop:Davies}. Let $I\subset \R$ be an open interval that contains $0$. Let $\th_0 \in [0,\pi/2)$. Then
\[
\Tr(\1_I^* \Pi^\h_{\b} \1_I)  \limt[\lvert \arg(\b)\rvert \leq \th_0]{|\b|}{\infty} 1.
\]
\end{corollary}
\begin{proof}
The reader who is not familiar with the trace of operators in infinite dimensional space may read, for instance, \cite[Chapter 10, \S1.3--1.4]{Kato95}.
We have $\Tr(\1_I^* \Pi^\h_{\b} \1_I) = \Tr( \1_I \1_I^* \Pi^\h_{\b})$.
Let $(\psi_k)_k$ be an orthonormal basis of $L^2(\R)$ such that $\psi_1 = \|\1_I\1_I^*\f_{\b,1}^\h\|^{-1}\1_I\1_I^*\f_{\b,1}^\h$. Then,
\begin{align*}
    \Tr(\1_I\1_1^*\Pi_\b^\h) 
    &= \sum_k \big\langle\1_I\1_I^*\Pi_\b^\h \psi_k, \psi_k\big\rangle\\
    &= \big\langle\1_I\1_I^*\Pi_\b^\h \psi_1, \psi_1\big\rangle\\
    &= \|\1_I\1_I^*\f_{\b,1}^\h\|^{-2}\frac {\big\langle \f_{\bar \b,1}^\h,\1_I\1_I^* \f_{\b,1}^\h \big\rangle} {\big\langle \f_{\bar \b,1}^\h , \f_{\b,1}^\h \big\rangle} \big\langle \f_{\b,1}^\h , \1_I\1_I^* \f_{\b,1}^\h\big \rangle\\
    &= \frac{\int_I \f_{\b,1}^\h(x)^2\diff x}{ \int_\R \f_{\b,1}^\h(x)^2 \diff x}.
\end{align*}
Since $0\in I$, the saddle point method proves that the right-hand side tends to $1$ as $|\beta|\to \infty$ and $\lvert\arg(\b)\rvert \leq \th_0$.
\end{proof}

 \section*{}

\subsection*{Acknowledgements}

We express our gratitude to our colleagues Joan Bellier-Mill\`es, Paulo Carrillo-Rouse and Joost Nuiten, for giving us the key arguments in algebraic topology to prove \cref{th-homotopy-2}.

This work has been partially supported by the ANR LabEx CIMI (under grant ANR-11-LABX-0040) within the French State Programme ``Investissements d'Avenir''.

\begin{microtypecontext}{expansion=bib}
% \printbibliography
\bibliographystyle{plain} % We choose the "plain" reference style
\bibliography{Grushin-gen}
\end{microtypecontext}
\end{document}

%% file: omega_neg_intro.tex
\begin{tikzpicture}[scale = 2.5]
  \draw (-2,0) rectangle (2,1);
% \draw[blue, fill=blue!50, thick] (-1.2,0.1) rectangle (2,0.3);
  \draw[fill = green!50] (-0.5,0) -- (-1.8,0) .. controls 
(-2,0.2) and (-2,0.6) .. (-1.8,0.7) .. controls 
(-1.6,0.9) and (-0.6,0.9) ..  (-0.6,0.7) .. controls 
(-0.6,0.6) and (-1.4,0.5) .. (-1.4,0.2) .. controls (-1.4,0.1) .. cycle;
  \draw[fill = green!50]  (1.6,0.4) .. controls 
(2,0.2) and (2,0.4) .. (1.8,1) --  (1,1) .. controls 
(1,0.3) and (1.5,0.4) .. cycle;
  \node (a) at (-0.5,1.2) {$\omega$};
  \draw[->] (a) -- (-1.3,0.7);
  \draw[->] (a) -- (1.6,0.7); 
  \draw[thick] (-1.2,0.2) -- (2,0.2) node[right]{$y_0$};
  \draw[dashed, thick, blue] (-1.2,0) node[below = -0.05]{$a$} -- (-1.2,0.2);
 \draw[->] (-2.2,0) -- (2.2,0) node[above]{$x$};
 \draw[->] (0,-0.2) -- (0,1.2) node[right]{$y$};
\end{tikzpicture}

%% file: omega-pos-intro.tex
\begin{tikzpicture}[scale=2.5]
\draw[fill=green!50] (-1.1,0) .. controls (-1.8,0.2) and (-1.5545,0.7659) .. (-1.3545,0.7659) .. controls (-1.1545,0.7659) and (-0.7992,0.3044) .. (0.0481,0.3126) .. controls (0.7,0.5) and (0.9,0.7) .. (-1.2,1) .. controls (0.5,1) and (1,1) .. (1,1) .. controls (1.2,0.9) and (0.9,0.2) .. (0.0444,0.192) .. controls (-0.6881,0.1878) and (-0.9,0.4) .. (-1.1451,0.5424) .. controls (-1.5,0.7) and (-1.4783,0.5111) .. (-1.3,0.4) .. controls (-0.6807,0.0438) and (0.9158,0.0626) .. (1,0) .. controls (-0.6,0) and (-1.1,0) .. (-1.1,0);
\draw[dashed, thick, blue] (-1.488,-0) node[below, black]{$\min_\T(\g_x)$} -- (-1.488,0.53);
\draw[dashed, thick, blue] (0.6,-0) node[below,black]{$\max_\T(\g_x)$} -- (0.6,0.585);

  \draw (-2,0) rectangle (2,1);
 \draw[->] (-2.2,0) -- (2.2,0) node[above]{$x$};
 \draw[->] (0,-0.2) -- (0,1.2) node[right]{$y$};
\draw (0.3,0.9) node[right]{$\omega$};
\draw[fill=white] (0.7,0.8) ellipse (0.05 and 0.05);
\draw[blue] (-0.5,0) .. controls (-1.4332,0.2846) and (-1.4225,0.299) .. (-1.4743,0.4381) .. controls (-1.5265,0.6283) and (-1.4,0.7) .. (-1.337,0.6876) .. controls (-1.2643,0.6719) and (-1.1059,0.5927) .. (-0.817,0.4359) .. controls (-0.5284,0.2964) and (-0.4264,0.2928) .. (-0.1871,0.2626) .. controls (0.0494,0.2351) and (0.2062,0.2708) .. (0.4758,0.4386) .. controls (0.7,0.6) and (0.7454,0.7164) .. (-0.5,1);
\draw[blue,->] (-0.5569,0.5726)node[above]{$\gamma$} -- (-0.7525,0.4054);
%\filldraw[opacity=0.75,pattern=dots,pattern color=red]  (-2,0) -- (-1.488,0) -- (-1.488,1) -- (-2,1) -- (-2,0);
%\filldraw[opacity=0.75,pattern=dots,pattern color=blue]  (0.6,0) -- (2,0) -- (2,1) -- (0.6,1) -- (0.6,0);

%\node[below, fill=white, fill opacity = 0.6, text opacity = 1] at (-1.8,0.5) {$\omega_-$};
%\node[below, fill=white, fill opacity = 0.6, text opacity = 1] at (1.4,0.5) {$\omega_+$};
\end{tikzpicture}

%% file: omega-corollary-intro.tex
\begin{tikzpicture}[scale=2.5]
 \draw (-2,0) rectangle (2,1);
 \draw[fill=green!50] (-1,1) .. controls (-1.9,0.7) and (-1.4,0.6) .. (-1.1,0.5) node (g1) {} .. controls (-0.8,0.4) and (0.5,0.3) .. (0.5,0.2) node (p) {} .. controls (0.5,0.1) and (-0.9,0.1) .. (-1,0) .. controls (-1,0) and (1.2,0) .. (1.2,0) .. controls (1.2,0.2) and (0.9,0.4) .. (0.6,0.5) node (g2) {} .. controls (0.3,0.6) and (-0.6,0.7) .. (-0.6,0.8) node(m){} .. controls (-0.6,0.9) and (1.1,0.9) .. (1.1,1) .. controls (1.1,1) and (-1,1) .. (-1,1);

 %\draw[thick, blue] (0,1) .. controls (-0.4,0.9) and (-0.7,0.9) .. (-0.7,0.8) .. controls (-0.7,0.7) and (-0.6,0.6) .. (-0.3,0.5) node (g) {} .. controls (0,0.4) and (0.6,0.3) .. (0.6,0.2) .. controls (0.6,0.1) and (0.4,0.1) .. (0,0);

 \node[fill = white, fill opacity = 0.4, text opacity = 1] at (-0.9352,0.7553) {$\omega$};
 
 %\draw[->, blue] (-1,0.3)node[left]{$\gamma$} -- (g.center);
 \draw[->] (-1.6,0.5)node[left]{$\gamma_1$} -- (g1.center);
 \draw[->] (1.3,0.6) node[right]{$\gamma_2$} -- (g2.center);

 \draw[dashed] (m.center) {} |- (0,0) node[pos=0.5,below]{$-\max_\T(\g_2^-)$};
 \draw[dashed] (p.center) {} |- (0,0) node[pos=0.5,below]{$\max_\T(\g_1^+)$};
%  \draw[dashed] (m.center) -| (0,0) node[pos = 0.5,right]{$y_-$};
%  \draw[dashed] (p.center) -| (0,0) node[pos = 0.5,left]{$y_+$};

 \draw[->] (-2.2,0) -- (2.2,0) node[above]{$x$};
 \draw[->] (0,-0.2) -- (0,1.2) node[right]{$y$};
\end{tikzpicture}

%% file: example3.tex
\begin{tikzpicture}[scale=2]
\draw (-1,0) -- (1,0) -- (1,1) -- (-1,1) -- (-1,0);
\draw[->] (-1.2,0) -- (1.2,0);
\draw[->] (0,-0.2) -- (0,1.2);
\draw[dashed, thick] (-0.5,0) -- (-0.5,0.375) ;
\draw[dashed, thick] (-0.75,0) -- (-0.75,0.625) ;
\draw  (-1,0) -- (-0.25,0) -- (-0.25,0.125) -- (-0.75,0.625) -- (-0.75,0.875) -- (-0.25,0.375) -- (-0.25,1) -- (-1,1) -- (-1,0);
\fill[opacity=0.5,color=gray,color=green]  (-1,0) -- (-0.25,0) -- (-0.25,0.125) -- (-0.75,0.625) -- (-0.75,0.875) -- (-0.25,0.375) -- (-0.25,1) -- (-1,1) -- (-1,0);
\path (-0.45,0.8) node {$\omega$};
\path (-0.5,-0.1) node {$-a$};
\path (-0.75,-0.1) node {$-b$};
%\path (-1.25,0.5) node {$c)$};
\end{tikzpicture}

%% file: example1.tex
\begin{tikzpicture}[scale=2]
\draw (-1,0) -- (1,0) -- (1,1) -- (-1,1) -- (-1,0);
\draw[->] (-1.2,0) -- (1.2,0);
\draw[->] (0,-0.2) -- (0,1.2);
\draw  (-0.75,0) -- (-0.75,0.75) -- (-0.25,0.75) -- (-0.25,0) -- (-0.75,0);
\draw  (0.75,1) -- (0.75,0.25) -- (0.25,0.25) -- (0.25,1) -- (0.75,1);
\fill[opacity=0.5,color=green]  (-0.75,0) -- (-0.75,0.75) -- (-0.25,0.75) -- (-0.25,0) -- (-0.75,0);
\fill[opacity=0.5,color=green]  (0.75,1) -- (0.75,0.25) -- (0.25,0.25) -- (0.25,1) -- (0.75,1);
\path (-0.5,0.25) node {$\omega$};
\path (0.5,0.75) node {$\omega$};
\path (-0.25,-0.1) node {$-a$};
%\path (-1.25,0.5) node {$a)$};
\end{tikzpicture}

%% file: example2.tex
\begin{tikzpicture}[scale=2]
\draw[->] (-1.2,0) -- (1.2,0);
\draw[->] (0,-0.2) -- (0,1.2);
\draw[dotted] (-0.5,0) -- (-0.5,0.5) ;
\draw[fill = green!50]  (-0.8,0) -- (-0.2,1) -- (-0.8,1) -- (-0.2,0) -- cycle;
\path (-0.6,0.15) node {$\omega$};
\path (-0.5,0.75) node {$\omega$};
\draw[thick, dashed] (-0.5,0.5) -- (-0.5,0) node[below] {$-a$};
\draw (-1,0) -- (1,0) -- (1,1) -- (-1,1) -- (-1,0);
%\path (-1.25,0.5) node {$b)$};
\end{tikzpicture}

%% file: omega0_left_right.tex
\begin{tikzpicture}[scale=2.5]
\draw[fill=green!50] (-1.1,0) .. controls (-1.8,0.2) and (-1.5545,0.7659) .. (-1.3545,0.7659) .. controls (-1.1545,0.7659) and (-0.7992,0.3044) .. (0.0481,0.3126) .. controls (0.7,0.5) and (0.9,0.7) .. (-1.2,1) .. controls (0.5,1) and (1,1) .. (1,1) .. controls (1.2,0.9) and (0.9,0.2) .. (0.0444,0.192) .. controls (-0.6881,0.1878) and (-0.9,0.4) .. (-1.1451,0.5424) .. controls (-1.5,0.7) and (-1.4783,0.5111) .. (-1.3,0.4) .. controls (-0.6807,0.0438) and (0.9158,0.0626) .. (1,0) .. controls (-0.6,0) and (-1.1,0) .. (-1.1,0);
\node[below] at (-1.488,-0) {$\min_\T(\g_x)$};
\node[below] at (0.6,-0) {$\max_\T(\g_x)$};
%\node[below] at (-1.43,-0) {$-a$};
%\node[below] at (1.43,-0) {$a$};

  \draw (-2,0) rectangle (2,1);
 \draw[->] (-2.2,0) -- (2.2,0) node[above]{$x$};
 \draw[->] (0,-0.2) -- (0,1.2) node[right]{$y$};
\draw (0.3,0.9) node[right]{$\omega$};
\draw[fill=white] (0.7,0.8) ellipse (0.05 and 0.05);
\draw[blue] (-0.5,0) .. controls (-1.4332,0.2846) and (-1.4225,0.299) .. (-1.4743,0.4381) .. controls (-1.5265,0.6283) and (-1.4,0.7) .. (-1.337,0.6876) .. controls (-1.2643,0.6719) and (-1.1059,0.5927) .. (-0.817,0.4359) .. controls (-0.5284,0.2964) and (-0.4264,0.2928) .. (-0.1871,0.2626) .. controls (0.0494,0.2351) and (0.2062,0.2708) .. (0.4758,0.4386) .. controls (0.7,0.6) and (0.7454,0.7164) .. (-0.5,1);
\draw[blue,->] (-0.5569,0.5726)node[above]{$\gamma$} -- (-0.7525,0.4054);
\filldraw[opacity=0.75,pattern=dots,pattern color=red]  (-2,0) -- (-1.488,0) -- (-1.488,1) -- (-2,1) -- (-2,0);
\filldraw[opacity=0.75,pattern=dots,pattern color=blue]  (0.6,0) -- (2,0) -- (2,1) -- (0.6,1) -- (0.6,0);

\node[below, fill=white, fill opacity = 0.6, text opacity = 1] at (-1.8,0.5) {$\omega_-$};
\node[below, fill=white, fill opacity = 0.6, text opacity = 1] at (1.4,0.5) {$\omega_+$};
\end{tikzpicture}

%% file: omega_neg.tex
\begin{tikzpicture}[scale = 2.5]
  \draw (-2,0) rectangle (2,1);
  \draw[blue, fill=blue!50, thick] (-1.2,0.1) rectangle (2,0.3);
  \draw[fill = green!50] (-0.5,0) -- (-1.8,0) .. controls 
(-2,0.2) and (-2,0.6) .. (-1.8,0.7) .. controls 
(-1.6,0.9) and (-0.6,0.9) ..  (-0.6,0.7) .. controls 
(-0.6,0.6) and (-1.4,0.5) .. (-1.4,0.2) .. controls (-1.4,0.1) .. cycle;
  \draw[fill = green!50]  (1.6,0.4) .. controls 
(2,0.2) and (2,0.4) .. (1.8,1) --  (1,1) .. controls 
(1,0.3) and (1.5,0.4) .. cycle;
  %\draw[fill=white] (-1.6,0.6) circle[radius=0.1];
  \node (a) at (-0.5,1.2) {$\omega$};
  \draw[->] (a) -- (-1.3,0.7);
  \draw[->] (a) -- (1.6,0.7); 
  \draw[thick] (-1.2,0.2) -- (2,0.2) node[right]{$y_0$};
%  \draw[thick, dashed] (0,0.2) -- (2,0.2) node[right]{$y_0$};
%   \draw[thick] (-1.4,0.2) -- (1.4,0.2) node[pos = 0.5, above right, fill = white, 
% fill opacity = 0.5, text opacity = 1, inner sep = 1pt]{$y_0$};
  %\draw[dashed, thick] (-1.4,0) node[below left = 0 and -0.2]{$-a$} -- (-1.4,0.2);
  %\draw[dashed, thick] (1.6,0) node[below = -0.02]{$b$} -- 
(1.6,0.2);
  \draw[dashed, thick, blue] (-1.2,0) node[below = -0.05]{$a$} -- (-1.2,0.2);
  %\draw[dashed, thick, blue] (1.4,0) node[below = -0.05]{$b$} -- (1.4,0.2);
 \draw[->] (-2.2,0) -- (2.2,0) node[above]{$x$};
 \draw[->] (0,-0.2) -- (0,1.2) node[right]{$y$};
\end{tikzpicture}

%% file: U.tex
\begin{tikzpicture}[scale=2.1]
\draw[fill=yellow!50]  (27:0.4) arc(27:45:0.4) -- (45:1) arc(45:387:1) -- cycle;
\node at(-0.4,0.4){$U$};
\draw[thick,->](-1.2,0) -- (1.2,0);
\draw[thick,->] (0,-1.2) -- (0,1.2);
\draw[very thick, |-|, blue] (27:1) arc (27:45:1) node[pos=0.5, right]{$W_0$};
%\draw[thick] (0,0) -- (36:1.3);
%\draw[thick, ->] (0:0.6) arc(0:36:0.6) node[pos=0.5, right]{$y_0$};
\draw[thin, gray] (0,0) circle[radius = 0.4];
\draw[<->] (0,0) -- (-110:0.4) node[pos = 0.5, left, 
fill=yellow!50, text=black, fill opacity = 0.5, inner sep 
= 0pt, text opacity = 1]{$\eu^{-a^2\!/2}$};
% \draw[pattern = dots, pattern color = red] (0,0) circle[radius=0.6];
% \draw[red,->] (-1.2,-0.8) node[left]{$D(0,e^{-T})$} -- (-0.2,-0.2);
% \draw[blue, thick, arrows={Rays[n=4]-}] (22:0.5) node[below]{$z_0$} -- 
% (22:1.4);
\end{tikzpicture}

%% file: runge-R.tex
\begin{tikzpicture}[scale=2]
\draw[fill=yellow!50]  (27:0.4) arc(27:45:0.4) -- (45:1) arc(45:387:1) -- cycle;
\node at(-0.6,0.6)[fill=white, fill opacity = 0.5, text opacity = 1, inner sep 
= 2pt]{$U$};
\draw[fill=red, fill opacity = 0.6] (0,0) circle[radius=0.6];
\draw[red!80!black, thick, ->] (-1,-0.9) node[below]{$D(0,\eu^{-T})$} -- (-0.2,-0.2);
\draw[blue, thick, arrows={Rays[n=4]-}] (36:0.5) node[below,xshift = 3,fill=white, 
fill opacity = 0.4, text opacity = 1, inner sep = 1pt]{$z_0$} -- 
(36:1.4);
\draw[thick,->](-1.2,0) -- (1.2,0);
\draw[thick,->] (0,-1.2) -- (0,1.2);
\end{tikzpicture}

%% file: D.tex
\begin{tikzpicture}[scale=2.7]
\draw[fill=yellow!50]  (27:0.4) arc(27:45:0.4) -- (45:1.1) arc(45:387:1.1) -- cycle;
\draw[fill=yellow!70!green]  (18:0.3) arc(18:54:0.3) -- (54:1) arc(54:378:1) -- cycle;
\draw[->] (-0.9,0.9) node[above]{$V$} -- (145:1.05);
\draw[gray, thin] (0,0) circle[radius = 0.3];
\draw[<->] (0,0) -- (-110:0.3);
\draw[thick,->](-1.2,0) -- (1.2,0);
\draw[thick,->] (0,-1.2) -- (0,1.2);
\node at(-0.3,0.45)[fill=yellow!80!green, fill opacity = 0.5, text opacity = 1, inner sep = 0.5pt]{$U$};
\path (-90:0.3) node[below, fill=yellow!80!green, text=black, fill opacity = 0.5, inner sep 
= 0pt, text opacity = 1]{$r=\eu^{-(1-\e)\dagmon(a)}$};
\end{tikzpicture}

%% file: runge.tex
\begin{tikzpicture}[scale=2]
\draw[fill=yellow!50]  (27:0.4) arc(27:45:0.4) -- (45:1.1) arc(45:387:1.1) -- cycle;
\node at(-0.6,0.6)[fill=white, fill opacity = 0.5, text opacity = 1, inner sep 
= 2pt]{$V$};
\draw[fill=red, fill opacity = 0.6] (0,0) circle[radius=0.6];
\draw[red!80!black, thick, ->] (-1,-0.9) node[below]{$D(0,\eu^{-q'(0)T(1+\e)})\qquad$} -- (-0.2,-0.2);
\draw[blue, thick, arrows={Rays[n=4]-}] (36:0.5) node[below,xshift = 3,fill=white, fill opacity = 0.4, text opacity = 1, inner sep = 1pt]{$z_0$} -- (36:1.4);
\draw[thick,->](-1.2,0) -- (1.2,0);
\draw[thick,->] (0,-1.2) -- (0,1.2);
\end{tikzpicture}

%% file: Delta_th.tex
\begin{tikzpicture}[baseline=(a),scale=1]
  \node (a) at (0,0) {};
  \path[fill=red!60] (75:3) -- (75:0.8)
    arc[radius=0.8, start angle= 75, end angle = -75] -- (-75:3)
    arc[radius=3, start angle=-75, end angle=75];
  \draw[red, thick] (75:3) -- (75:0.8)
    arc[radius=0.8, start angle= 75, end angle = -75] -- (-75:3);
  \node at (-65:1.25) [right]{$\sector[\th]\setminus D(0,r_\th)$};
  \draw[->](-1,0) -- (3,0);
  \draw[->](0,-2.7) -- (0,2.7);
  \draw[thick] (0.8,0.2) -- +(0,-0.4) 
    node[below,fill=white, text=black, fill opacity = 0.5, inner sep = 1pt, text opacity = 1] {$r_\theta$};;
  \draw[->, thick] (1.2,0) arc[start angle = 0, end angle = 75, radius = 1.2];
  \node at (30:1.2)[above right] {$\theta$};
\end{tikzpicture}

%% file: omega_corollary.tex
\begin{tikzpicture}[scale=2.5]
 \draw (-2,0) rectangle (2,1);
 \draw[fill=green!50] (-1,1) .. controls (-1.9,0.7) and (-1.4,0.6) .. (-1.1,0.5) node (g1) {} .. controls (-0.8,0.4) and (0.5,0.3) .. (0.5,0.2) node (p) {} .. controls (0.5,0.1) and (-0.9,0.1) .. (-1,0) .. controls (-1,0) and (1.2,0) .. (1.2,0) .. controls (1.2,0.2) and (0.9,0.4) .. (0.6,0.5) node (g2) {} .. controls (0.3,0.6) and (-0.6,0.7) .. (-0.6,0.8) node(m){} .. controls (-0.6,0.9) and (1.1,0.9) .. (1.1,1) .. controls (1.1,1) and (-1,1) .. (-1,1);

 \draw[thick, blue] (0,1) .. controls (-0.4,0.9) and (-0.7,0.9) .. (-0.7,0.8) .. controls (-0.7,0.7) and (-0.6,0.6) .. (-0.3,0.5) node (g) {} .. controls (0,0.4) and (0.6,0.3) .. (0.6,0.2) .. controls (0.6,0.1) and (0.4,0.1) .. (0,0);

 \node[fill = white, fill opacity = 0.4, text opacity = 1] at (-0.9352,0.8553) {$\omega$};
 
 \draw[->, blue] (-1,0.3)node[left]{$\gamma$} -- (g.center);
 \draw[->] (-1.6,0.5)node[left]{$\gamma_1$} -- (g1.center);
 \draw[->] (1.3,0.6) node[right]{$\gamma_2$} -- (g2.center);

 \draw[dashed] (m.center) {} |- (0,0) node[pos=0.5,below]{$a_-$};
 \draw[dashed] (p.center) {} |- (0,0) node[pos=0.5,below]{$a_+$};
 \draw[dashed] (m.center) -| (0,0) node[pos = 0.5,right]{$y_-$};
 \draw[dashed] (p.center) -| (0,0) node[pos = 0.5,left]{$y_+$};

 \draw[->] (-2.2,0) -- (2.2,0) node[above]{$x$};
 \draw[->] (0,-0.2) -- (0,1.2) node[right]{$y$};
\end{tikzpicture}

%% file: covering.tex
\begin{tikzpicture}[scale=1.5]
\begin{scope}
 \draw (-2,0) rectangle (2,1);
 \draw[->] (-2.2,0) -- (2.2,0) node[above]{$x$};
 \draw[->] (0,-0.2) -- (0,1.2) node[right]{$y$};
 \draw[thick, blue, postaction={decorate, 
      decoration={markings, mark = between positions 0.2 and 0.9 step 0.28 with {\arrow{>}}}}] 
  (-2,0.5) node(ya){} to[in=-60, out = 40] node[pos= 0.85](n1){} (0.5,1)  
  (0.5,0) to[in=120, out=120] node(n2)[pos = 0.6]{} (2,0.25) node(yb){};
  
  \draw (1.2,1.4) node[text = blue](g){$c_1$};
  \draw[thick, blue, ->] (g) -- (n1.center);
  \draw[thick, blue, ->] (g) -- (n2.center);
 \draw[thick, dashed] let \p1 = (ya) in (ya.center) -- (0,\y1) node[right]{$y_0$};
 \draw[thick, dashed] let \p1 = (yb) in (yb.center) -- (0,\y1) node[left]{$y_1$};
 \draw[->] (-2.2,0) -- (2.2,0) node[above]{$x$};
 \draw[->] (0,-0.2) -- (0,1.2) node[right]{$y$};
\end{scope}
\begin{scope}[shift = {(5,-0.7)}]
 \draw[thick, blue, postaction={decorate, 
      decoration={markings, mark = between positions 0.1 and 0.9 step 0.15 with {\arrow{>}}}}] 
  (-2,0.5) node(tya){} to[in=-60, out = 40] (0.5,1)  to[in=120, out=120] (2,1.25) node(tyb){}
  -- (2,2.25) to[in=120, out=120] (0.5,2) to[out=-60, in= 40] node(n){} (-2,1.5)
  --cycle;
  \draw (-1.2,2.1) node[text = blue](gg){$c_2$};
  \draw[thick, blue, ->] (gg) -- (n.center);
  
 \draw[thick, dashed] let \p1 = (tya) in (tya.center) -- (0,\y1) node[right]{$\tilde y_0$};
 \draw[thick, dashed] let \p1 = (tyb) in (tyb.center) -- (0,\y1) node[left]{$\tilde y_1$};
 \draw[->] (-2.2,0) -- (2.2,0) node[above]{$x$};
 \draw[->] (0,-0.2) -- (0,2.6) node[right]{$y$};
\end{scope}
\end{tikzpicture}